\documentclass[11pt,leqno]{amsart}
\usepackage{times,amsmath,amsthm,amssymb,latexsym,amscd}
\usepackage[all]{xy}
\usepackage{epsfig}

\CompileMatrices
\newtheoremstyle{fancy}{}{}{\itshape}{}{\textsc\bgroup}{.\egroup}{ }{}
\newtheoremstyle{fanci}{}{}{\rm}{}{\textsc\bgroup}{.\egroup}{ }{}
\theoremstyle{fancy}
\newcounter{intro}

\numberwithin{equation}{section} \swapnumbers
\newtheorem{cor}[equation]{Corollary}
\newtheorem{lem}[equation]{Lemma}
\newtheorem{prop}[equation]{Proposition}
\newtheorem{thm}[equation]{Theorem}
\theoremstyle{fanci}
\newtheorem{dfn}[equation]{Definition}
\newtheorem{exa}[equation]{Example}

\newtheorem{rem}[equation]{Remark}
\newtheorem{rems}[equation]{Remarks}
\newtheorem{onb}[equation]{Orthonormal frame bundle}
\newtheorem{nota}[equation]{Notation}
\newtheorem{notarem}[equation]{Notation and Remarks}
\newtheorem{bk}[equation]{Computation of the $b_k$}
\newtheorem{specpf}[equation]{Proof of Theorem \ref{thm.trace} in special case}
\newcommand{\cref}[1]{Corollary~\ref{#1}}

  \def\R{\mathbb R} \def\Z{\mathbb{Z}}

\newcommand{\spec}{\operatorname{spec}}
\newcommand{\supp}{\operatorname{supp}}

\newcommand{\vol}{\operatorname{vol}}

\def\O{{\mathcal{O}}}
\def\a{{\alpha}}

\def\g{{\gamma}}
\def\ga{{G_\a}}
\def\hm{{H^{(m)}}}
\def\tham{{\widetilde{H}^{(m)}_\a}}
\def\ham{{H^{(m)}_\a}}
\def\tha{{\widetilde{H}_\a}}
\def\ha{{H_\a}}

\def\phia{{\eta_\a}}
\def\psia{{\psi_\a}}
\def\too{{t\to 0^+}}
\def\dx{{\delta_x}}
\def\ua{{U_\a}}
\def\va{{V_\a}}
\def\wa{{W_\a}}
\def\wt{\widetilde}
\def\tu{{\widetilde{U}}}
\def\tul{{\widetilde{u}}}
\def\tua{{\widetilde{U}_\a}}
\def\tva{{\widetilde{V}_\a}}
\def\twa{{\widetilde{W}_\a}}

\def\tpsia{{\widetilde{\psi}_\a}}

\def\tx{{\widetilde{x}}}
\def\ty{{\widetilde{y}}}
\def\tn{{\widetilde{N}}}
\def\heat{{\left(\frac{\partial}{\partial t}+\Delta_x\right)}}
\def\theat{{\left(\frac{\partial}{\partial t}+\widetilde{\Delta}_\tx\right)}}

\def\tfa{{\widetilde{f}_\a}}

\def\bs{{\backslash}}

\def\pa{{p_\a}}
\def\tpa{{\widetilde{p}_\a}}

\def\shag{{(4\pi t)^{-n/2}e^{-d(\tx,\g(\ty))^2/4t}}}

\def\shagx{{(4\pi t)^{-n/2}e^{-d(\tx,\g(\tx))^2/4t}}}

\def\sha{{(4\pi t)^{-n/2}e^{-d(\tx,\ty)^2/4t}}}
\def\eha{{(u_0(\tx,\ty)+\dots +t^mu_m(\tx,\ty))}}

\def\ehagi{{(u_0(\tx,\g(\ty))+ tu_1(\tx,\g(\ty))+\dots)}}
\def\ehagxi{{(u_0(\tx,\g(\tx))+ tu_1(\tx,\g(\tx))+\dots)}}
\def\ehaxi{{(u_0(\tx,\tx)+ tu_1(\tx,\tx)+\dots)}}
\def\gn{{\Iso(N)}}
\def\e{{\epsilon}}

\def\heatg{{\shagx\ehagxi}}

\def\om{{\Omega}}
\def\tn{{\widetilde{N}}}

\def\pa{{p_\a}}

\def\gns{{\Iso^{\operatorname{max}}(N)}}

\def\rp{{\R_+}}
\def\rps{{\R_+^*}}

\def\supp{{\text{supp}}}
\def\tp{{\tilde{p}}}

\def\tv{{\widetilde{V}}}
\def\tetaa{{\widetilde{\eta}_\a}}
\def\etaa{{\eta_\a}}
\def\teta{{\widetilde{\eta}}}
\def\tb{{\widetilde{b}}}
\def\ct{{c}}

\def\Iso{{\operatorname{{Iso}}}}
\def\ismax{{{\Iso}^{\operatorname{max}}}}
\def\Fix{{\operatorname{{Fix}}}}

\begin{document}

\title[Asymptotic expansion of the heat kernel for orbifolds]{Asymptotic
expansion of the heat kernel for orbifolds}

\author[Dryden]{Emily B. Dryden}
\address[Dryden]{Department of Mathematics, Bucknell University, Lewisburg, PA 17837}
\email{ed012@bucknell.edu}

\author[Gordon]{Carolyn S. Gordon}
\address[Gordon]{Department of Mathematics, Dartmouth College, Hanover, 
New Hampshire 03755}
\email{csgordon@dartmouth.edu}

\author[Greenwald]{Sarah J. Greenwald}
\address[Greenwald]{Department of  Mathematics, Appalachian 
State University, Boone, NC 28608}
\email{greenwaldsj@appstate.edu}

\author[Webb]{David L. Webb}
\address[Webb]{Department of Mathematics, Dartmouth College, Hanover, 
New Hampshire 03755}
\email{david.l.webb@dartmouth.edu}

\thanks{Gordon and Webb were
supported in part by NSF grants DMS-0072534, DMS-0306752, and DMS-0605247; 
Greenwald was supported in part by NSF ROA grants 0072533 and
9972304.}

\begin{abstract}
We study the relationship between the geometry and the Laplace spectrum of a
Riemannian orbifold $\O$ via its heat kernel; as in the manifold case, the
time-zero asymptotic expansion of the heat kernel furnishes geometric
information about $\O$.  In the case of a \emph{good} Riemannian orbifold (i.e.,
an orbifold arising as the orbit space of a manifold under the action of a
discrete group of isometries), H. Donnelly \cite{D1} proved the existence of the
heat kernel and constructed the asymptotic expansion for the heat trace.  We
extend Donnelly's work to the case of general compact orbifolds.  Moreover, in
both the good case and the general case, we express the heat invariants in a
form that clarifies the asymptotic contribution of each part of the singular set
of the
orbifold.  We calculate several terms in the asymptotic expansion explicitly
in the case of two-dimensional orbifolds; we use these terms to prove that
the spectrum distinguishes elements within various classes of two-dimensional
orbifolds.

\end{abstract}

\maketitle

\tableofcontents

\section{Introduction}

Roughly speaking, a topological \emph{orbifold} is a space locally homeomorphic to
an orbit space of a finite group action on $\R^n$.  A smooth orbifold consists of a Hausdorff second countable
topological space together with an atlas of
coordinate charts realizing such local homeomorphisms and satisfying 
compatibility conditions (see Section
\ref{sec.oflds}). Orbifolds
were introduced by Satake, then studied by Thurston because of their utility
in the investigation of three-manifolds (e.g., a Seifert fibred three-manifold
is naturally a generalized circle bundle over a two-orbifold); today, orbifolds
arise naturally in diverse branches of mathematics and physics, including
symplectic
geometry, string theory, and vertex operator algebras.

We will be interested in orbifolds from a spectral-theoretic point of view.  
An orbifold endowed with a metric structure
is a \emph{Riemannian orbifold}.  As in the manifold case, associated with every
Riemannian metric is a Laplace operator acting on smooth functions on the
orbifold.  In the case of closed orbifolds, the Laplacian has a discrete
spectrum.
We study the relationship between the geometry and the Laplace spectrum of a
closed
orbifold via its heat kernel; as in the manifold case, the time-zero asymptotic
expansion of the heat kernel furnishes geometric information about the orbifold.

Orbifolds began appearing sporadically in the spectral theory literature in the
early 1990s, and have received more concentrated attention in the last
five years.  C. Farsi \cite{Far01} showed that the spectrum of an orbifold
determines its volume by proving that Weyl's asymptotic formula holds for
orbifolds. Dryden and A. Strohmaier \cite{EBDStroh} showed that for a compact,
negatively curved
two-dimensional orbifold, the Laplace spectrum determines both the length
spectrum
and the orders of the singular points and vice versa;
on the other hand, P. Doyle and J.P. Rossetti \cite{DoyRo} gave (disconnected) 
examples of isospectral flat two-dimensional orbifolds with different length
spectra 
and orders of singular points. 
Further investigations of the relationship between the lengths
of closed geodesics and the spectrum were carried out by E. Stanhope and A.
Uribe in \cite{SU}.  It is natural to ask about the singularities that can
appear in an
isospectral family of orbifolds.  Stanhope \cite{Sta05} showed that, in
general, there can be at most finitely many isotropy types (up to isomorphism)
in a set of isospectral Riemannian orbifolds that share a uniform lower bound on
Ricci curvature.  On the other hand, N. Shams, Stanhope, and Webb \cite {SSW}
constructed arbitrarily large (finite) isospectral sets of orbifolds
satisfying this curvature condition whose isotropy types differ.  Rossetti, D. Schueth, and M. Weilandt 
\cite{Sch} recently constructed a pair of isospectral Riemannian orbifolds whose
isotropy types have different orders.  

For Riemannian manifolds, the asymptotic expansion of the heat kernel can be
used to relate the geometry of the manifold to its spectrum.  From the so-called
\emph{heat invariants} appearing in the asymptotic expansion, one can tell the
dimension, the volume, and various quantities involving the curvature of the
manifold.  The heat kernel has been studied in various analogous or more general
settings (e.g. \cite{BG90,BL96,BS87,D2,Gil,Ric}).  

In the case of a \emph{good} Riemannian orbifold (i.e., an orbifold arising as
the orbit space of a manifold under the action of a discrete group of
isometries), H. Donnelly
\cite{D1} proved the existence of the heat kernel and constructed the asymptotic
expansion for the heat trace.  We extend Donnelly's work to the case of general
compact orbifolds.   Moreover, in both the good case and the general case, we 
express the heat invariants in a form that clarifies
the asymptotic contribution of each part of the singular set of the orbifold.   

We calculate several terms in the asymptotic expansion explicitly
in the case of two-dimensional orbifolds; we use these terms to prove that the
spectrum
distinguishes elements within various classes of two-dimensional orbifolds.  In
particular, within the class of all two-dimensional orbifolds with nonnegative
Euler characteristic, the spectrum is a complete topological invariant.  
Additional results are obtained for triangular pillow orbifolds endowed with a
hyperbolic
structure, and for nonorientable two-dimensional orbifolds.  

The paper is organized as follows.  In Section \ref{sec.oflds}, we give the
background necessary for the rest of the paper, recalling several
results that clarify the structure of the singular locus of an orbifold.
Section \ref{sec.heatkernel} is devoted to the construction of the heat kernel
on an arbitrary closed Riemannian orbifold by means of the construction of a
parametrix.  The existence of the heat kernel for closed orbifolds was shown previously by Y.-J. Chiang \cite{Chi90}; existence also follows from more general
results for the heat kernel on Riemannian foliations (see \cite{Ric}).  However,
we give a different construction in order to express the heat kernel in
a convenient form that will allow us, in Section \ref{sec.asymptotics}, to
generalize Donnelly's asymptotic expansion.  Section \ref{sec.apps} is devoted
to various applications of the heat expansion, including those mentioned in the
previous paragraph. \\

\subsection*{Acknowledgments}  The initial inspiration for this project came from Shunhui Zhu, and the authors thank him for sharing his curiosity.  We also thank Iosif Polterovich for helpful
discussions and encouragement, and Alejandro Uribe and Liz Stanhope for alerting us to the frame bundle approach. 
The first named author acknowledges the Centre Interfacultaire Bernoulli in Lausanne, Switzerland, and the Centro de An\'{a}lise Matem\'{a}tica, Geometria e Sistemas Din\^{a}micos, Instituto Superior T\'{e}cnico, in Lisbon, Portugal, for their support during her work on this project.

%%%%%%%%%%%%%%%%%%%%%%%%%%%%%%%%%%%%%%%%%%%%%%%%%%%%%%%%%%%%
%%%%%%%%%%%%%%%%%%%%%%%%%%%%%%%%%%%%%%%%%%%%%%%%%%%%%%%%%%%%
\section{Orbifolds and their singular sets}\label{sec.oflds}
%%%%%%%%%%%%%%%%%%%%%%%%%%%%%%%%%%%%%%%%%%%%%%%%%%%%%%%%%%%%
%%%%%%%%%%%%%%%%%%%%%%%%%%%%%%%%%%%%%%%%%%%%%%%%%%%%%%%%%%%%
  
\begin{dfn}\label{def.orb} \text{   }
\begin{enumerate}
\item An \emph{orbifold chart} on a topological space $X$ consists of a
connected open subset $\widetilde{U}$ of $\R^n$, a finite group $G_U$ acting on 
$\widetilde{U}$ by diffeomorphisms, and a mapping $\pi_U$ from $\widetilde{U}$
onto an open subset $U$ of $X$ inducing a homeomorphism from the orbit space
$G_U\backslash\widetilde{U}$ onto $U$.  We will always assume that the group
$G_U$ acts effectively on $\widetilde{U}$. 

\item An \emph{embedding} $\lambda:(\widetilde{U},G_U,\pi_U)\to
(\widetilde{U}', G_{U'},\pi_{U'})$ between orbifold charts with 
$U\subseteq U'$ is a smooth embedding
$\lambda:\widetilde{U}\to\widetilde{U}'$ such that 
the diagram
$$\xymatrix{
{\wt{U}\,}\ar^{\lambda}@{>->}[r]\ar_{\pi_U}@{->>}[d]&
{\wt{U'}}\ar^{\pi_{U'}}@{->>}[d]\\
{G_U\bs\wt{U}}\ar_{\sim}@{->}[d]&
{G_{U'}\bs\wt{U'}}\ar^{\sim}@{->}[d]\\
U\ar@{^{(}->}[r]&{U'}
}$$
commutes.
Two charts $(\widetilde{U},G_U,\pi_U)$ and
$(\widetilde{U}', G_{U'},\pi_{U'})$
on $X$ are said to be \emph{compatible} if, for each point $x\in U\cap U'$, there
exists an orbifold chart 
$(\widetilde{U}'', G_{U''},\pi_{U''})$ with $U''\subset U\cap U'$
and smooth embeddings 
$(\widetilde{U}'', G_{U''},\pi_{U''}) \to
(\widetilde{U},G_U,\pi_{U})$ and $(\widetilde{U}'', G_{U''},\pi_{U''})\to
(\widetilde{U}', G_{U'},\pi_{U'})$.  

\item An $n$-dimensional \emph{orbifold atlas}
$\mathcal A$ on $X$ is a compatible family of $n$-dimensional orbifold charts
whose images form a covering of $X$.  A \emph{refinement} $\mathcal{A}'$ of an
orbifold atlas $\mathcal{A}$ is an orbifold atlas each of whose
charts embeds into a chart of $\mathcal{A}$.  Two orbifold atlases are said to
be \emph{equivalent} if they have a common refinement.  Every orbifold atlas is
equivalent to a unique maximal one.  An \emph{orbifold} is a Hausdorff, second
countable topological space together with a maximal orbifold atlas.

\item Let $\O$ be an orbifold.  A point $x$ of
$\O$ is said to be \emph{singular} if for some (hence every) orbifold chart
$(\widetilde{U},G_U,\pi_U)$ about $x$, the points in the inverse image of $x$ in
$\widetilde{U}$ have nontrivial isotropy in $G_U$.  The isomorphism class of
the isotropy group, called the \emph{abstract isotropy type} of $x$, 
 is independent both of the choice of point in the inverse
image of $x$ in $\tu$ and of the choice of chart $(\widetilde{U},G_U,\pi_U)$
about $x$.   
Points that are not singular are called \emph{regular}.  
\end{enumerate}
\end{dfn}
  
\begin{rems}\label{rem.embed} \text{   }
\begin{enumerate}

\item The notion of orbifold generalizes slightly the notion of
$V$-manifold introduced by Satake \cite{Sat57}; $V$-manifolds are
orbifolds for
which the singular set has codimension at least two.

\item Given an embedding $\lambda:(\widetilde{U},G_U,\pi_U)\to
(\widetilde{U}', G_{U'},\pi_{U'})$ as in Definition \ref{def.orb}, there exists
a homomorphism $\tau:G_U\to G_{U'}$ such that 
$\lambda\circ \g=\tau(\g)\circ\lambda$ for
all $\g\in G_U$.  This was proven by Satake for $V$-manifolds and was
generalized to orbifolds by Moerdijk and Pronk \cite{MP97}.  From our
convention that the group actions on each chart are effective, it follows that
the homomorphism $\tau$ is injective.

\item  There are some subtle differences among the definitions of
\emph{orbifold} in the literature; in particular, some authors do not require that the
finite group actions in the orbifold charts be effective.  The definition we
use is that in \cite{MP97}. 
While these distinctions become significant when formulating the correct notion
of the category of smooth orbifolds, they play no role in our computations.
\end{enumerate}
\end{rems}
  
An orbifold is said to be \emph{good} if it is the orbit space of a 
manifold under the smooth action of a discrete group; otherwise, it is said to
be \emph{bad}.  In particular, every point in an orbifold has a 
neighborhood that is a good orbifold.  Even a bad orbifold can be expressed
globally as the quotient of a manifold by a group action, although not a
discrete group action.  This is done by introducing a Riemannian structure and
constructing the ``bundle'' of orthonormal frames.  The orthonormal frame bundle
is actually an
\emph{orbibundle}, the appropriate notion of vector bundle over an orbifold as
described in \cite{M02} or \cite{Thu80}.   
However, its total space is a smooth manifold with an action of the orthogonal
group, and one recovers the original orbifold as the orbit space of this
orthogonal action.  We will review this construction below after first
establishing some notation in the setting of arbitrary group actions.

\begin{dfn}\label{def.orbittype}  (See \cite[Chapter 2]{DK}.)
\begin{enumerate}

\item Consider a smooth proper action of a Lie group $H$ on a smooth
manifold $M$.  For $\tx\in M$, let ${\Iso}_H(\tx)$ denote the subgroup of $H$ that fixes $\tx$.   
Define an equivalence relation on $M$ by $\tx\equiv \ty$ if ${\Iso}_H(\tx)$ and
${\Iso}_H(\ty)$ are conjugate.  Each equivalence class is called an
$H$-\emph{orbit type}.  Note that the equivalence classes are invariant under
the action of $H$. 

We will say that the $H$-orbit type of $\tx$ \emph{dominates} that of $\ty$ if
${\Iso}_H(\tx)$ is conjugate to a subgroup of ${\Iso}_H(\ty)$.

\item Let $\pi:M\to H\bs M$ be the projection onto the orbit space.  
Let $p\in H\bs M$.  As $\tp$ ranges over the $H$-orbit $\pi^{-1}(p)$ in 
$M$, the stabilizer ${\Iso}_H(\tp)$ ranges over a conjugacy class of
subgroups of $H$.  We will denote this conjugacy class of subgroups by 
${\Iso}_H(p)$ and refer to it as the $H$-\emph{isotropy type} of $p$.  Define an
equivalence relation on $H\bs M$ by  $p\equiv q$ if
${\Iso}_H(p)={\Iso}_H(q)$.  The equivalence classes will be called
$H$-\emph{isotropy equivalence classes}.   Note that $\pi$ carries points of the
same 
$H$-orbit type in $M$ to points of the same
$H$-isotropy equivalence class in $H\bs M$.  

We will say that the $H$-isotropy equivalence class of $p$ \emph{dominates} that
of $q$ if
the groups making up the conjugacy class ${\Iso}_H(p)$ are conjugate to
subgroups of those in ${\Iso}_H(q)$.

By an abuse of notation, we will write $|\Iso_H(p)|$ to mean the order of each
of the groups making up the isotropy type $\Iso_H(p)$.  We will refer to this
quantity as the order of the $H$-isotropy at $p$.
\end{enumerate}
\end{dfn}

\begin{dfn}\label{def.riem}  \text{   }
\begin{enumerate}

\item A \emph{Riemannian structure} on an orbifold $\O$
is an assignment to each orbifold chart $(\widetilde{U},G_U,\pi_U)$ 
of a $G_U$-invariant Riemannian metric $g_{\widetilde{U}}$ on $\wt{U}$
satisfying the compatibility condition that each embedding $\lambda$ appearing
in
Definition \ref{def.orb} is isometric.  Every orbifold admits Riemannian
structures.  

\item We will say that an orbifold chart $(\tu,G_U,\pi_U)$ on a
Riemannian orbifold $\O$ is a \emph{distinguished} chart of radius $r$ if $\tu$
is a convex geodesic ball of radius $r$.   In this case,  $U$ is a convex
geodesic ball in $\O$. The entire group $G_U$ fixes the center 
$\tp$ of $\tu$, so the abstract isotropy type of $p:=\pi_U(\tp)$ is represented
by $G_U$.  

\end{enumerate}
\end{dfn}

\begin{rem}\label{convex.rad}  Recall that for a Riemannian manifold $M$ and a
point $p\in M$, the \emph{convexity radius at $p$} is the largest positive real
number $r(p)$ for which the geodesic ball of radius $\epsilon$ about $p$ is
geodesically convex for all $\epsilon <r(p)$.  If $M$ is compact, the
infimum $r$ of $\{r(p):p\in M\}$ is positive and is called the \emph{convexity
radius} of $M$.  For a point $p$ in an  orbifold $\O$, we may define the
convexity radius at $p$ to be the largest real number $r(p)$ such that $\O$
admits a distinguished chart of radius $\epsilon$ centered at $p$ for all
$\epsilon <r(p)$.  It is immediate that $r(p)$ is positive.  
Moreover, if $\O$ is compact, then an elementary argument shows that the infimum
$r$ of $\{r(p):p\in\O\}$ is positive; $r$ is called the \emph{convexity radius}
of $\O$.
\end{rem}
 
\begin{onb}\label{onb} We give a brief description of the orthonormal frame
bundle of a Riemannian  orbifold.  See \cite{AR03} for more details.  First
consider a
good Riemannian orbifold $\O=G\bs M$, where $M$ is a Riemannian manifold and $G$
is a discrete group acting by
isometries on  $M$. 
Let $F(M)\to M$ be the orthonormal frame bundle of  $M$. 
Each element $\g\in G$, being an isometry of $M$, induces a diffeomorphism
$\g_*$ of
$F(M)$ carrying fibers to fibers; thus we obtain an action of $G$ on $F(M)$
covering the action of $G$ on $M$.  
The orthonormal frame bundle $F(\O)$ of $\O$ is defined to be 
$G\bs F(M)\to\O$.  The \emph{fiber} of $F(\O)\to\O$ over $x\in\O$ is the 
preimage of $x$ in $G\bs F(M)$.  The right action of $O(n)$ on $F(M)$ commutes
with the left action of $G$, and hence descends to a right $O(n)$-action on
$F(\O)$.
For a bad orbifold, the orthonormal frame bundle
is defined in such a way that its restriction to any good neighborhood 
$U\cong G_U\bs\widetilde{U}$ is the orthonormal frame bundle of the good
orbifold $U$.
  
The orthonormal frame bundle of $\O$ is a smooth manifold as well as an
orbibundle on which the orthogonal group $O(n)$ acts smoothly on the right,
preserving
fibers.  In particular, the orbifold $\O$  is the orbit space $F(\O)/O(n)$ of
the right action of
$O(n)$ on the manifold $F(\O)$.   
\end{onb}

\begin{notarem}\label{localisot} Let $\O$ be an orbifold.  Endow $\O$ with a
Riemannian metric and let $F(\O)$ be the associated orthonormal frame bundle as
in \ref{onb}.  The $O(n)$-action on the fiber of $F(\O)$ over a point $x\in\O$
is free if and only if
$x$ is a regular point of $\O$.  In particular, for each singular point $x$ of
$\O$, viewed as an element of $F(\O)/O(n)$, the $O(n)$-isotropy type of $x$ is a
non-trivial conjugacy class $\Iso_{O(n)}(x)$ of subgroups of $O(n)$. (See
Definition \ref{def.orbittype}.)   

Our realization of the orbifold $\O$ as a global quotient of a manifold (namely
$F(\O)$) by an action of $O(n)$ depends, of course, on the Riemannian metric. 
However, it is not difficult to show that the conjugacy class $\Iso_{O(n)}(x)$
of subgroups of $O(n)$ is actually independent of the choice of Riemannian
metric used in the construction.    This conjugacy class of subgroups of $O(n)$
will henceforth be denoted $\Iso(x)$ (without the subscript $O(n)$ except when
needed for clarity) and will be referred to as the \emph{isotropy type} of the
singular point $x$ of $\O$.  Its cardinality $|\Iso(x)|$ will be called the
\emph{order} of the isotropy at $x$.  Similarly, the equivalence classes of elements of
$\O$ with the same isotropy type will be called \emph{isotropy equivalence classes},
without mention of $O(n)$.

The subgroups  of $O(n)$ in the conjugacy class $\Iso(x)$ lie in the isomorphism
class defined by the  abstract isotropy type of $x$ given in Definition
\ref{def.orb}.   Indeed, let  $(\tu,G_U,\pi_U)$ be an orbifold chart with $x\in
U$.  Let $\tx\in \tu$
with $\pi_U(\tx)=x$.  With respect to any choice of Riemannian metric on $\O$
(and associated Riemannian metric on $\tu$), the group $G_U$ acts isometrically
on  $\tu$ and thus acts on the  left
on the orthonormal frame bundle $F(\tu)$.  The subgroup ${\Iso}_{G_U}(\tx)$
leaves invariant the fiber of $F(\tu)$ over $\tx$.  For each $q$ in this fiber,
define a homomorphism $\sigma_q:{\Iso}_{G_U}(\tx)\to O(n)$ by the condition
$\gamma(q) =(q)\sigma_q(\gamma)$ where $\gamma(\cdot)$ and 
$(\cdot)\sigma_q(\gamma)$ denote the left action of $\gamma$ and right action of
$\sigma_q(\gamma)\in O(n)$ on the fiber.  By \ref{onb}, the restriction of
$F(\O)$ to $U$ is given by $G_U\bs F(\tu)$.  Letting $\rho:F(\tu)\to G_U\bs
F(\tu)$ be the projection, then $\sigma_q$ maps  ${\Iso}_{G_U}(\tx)$
isomorphically to the stabilizer of $\rho(q)$ in $O(n)$.  This stabilizer is a
representative of the conjugacy class ${\Iso}(x)$, while ${\Iso}_{G_U}(\tx)$
represents the abstract isotropy type of $x$.\end{notarem}

\begin{dfn} \label{def.strat} A smooth \emph{stratification} of a manifold or
orbifold $M$ is a locally finite partition of $M$ into locally closed
submanifolds, called the \emph{strata}, satisfying the following condition:  For
each stratum $N$, the closure of $N$ is the union of $N$ with a collection of
lower dimensional strata.   
\end{dfn}

\begin{rems}\label{rem.strat} \text{ }
\begin{enumerate}
\item For any stratification of an orbifold (or manifold) $\O$, the strata of
maximal dimension are open in $\O$ and their union has full measure in $\O$.
\item The stratifications that we will discuss below are Whitney
stratifications. 
As we will not explicitly use the additional properties of Whitney
stratifications here, we omit the definition and refer the reader to \cite{DK}. 
The notion of Whitney stratification can be defined in the more general setting
of spaces that can be at least locally embedded in a smooth manifold.  As
discussed in \cite{DK}, the orbit space of a proper Lie group action on a smooth
manifold has this property.
\end{enumerate}
\end{rems}

\begin{prop}\cite{DK}\label{prop.strat}  
Given a smooth action of a Lie group $H$ on a manifold $M$, then:
\begin{enumerate}
\item The connected components of the $H$-orbit types form a Whitney
stratification
of $M$.  The closure of a stratum $\tn$ is made up of the union of $\tn$ with a
collection of lower dimensional strata, each lying in an $H$-orbit type strictly
dominated by that of $\tn$.

\item The connected components of the $H$-isotropy equivalence classes in $H\bs
M$ form a
Whitney stratification of $H\bs M$.  The closure of a stratum $N$ is made up of
the union of $N$ with a collection of lower dimensional strata, each lying in an
$H$-isotropy equivalence class strictly dominated by that of $N$.
The map $\pi$ carries each stratum in $M$ onto a stratum in $H\bs M$.

\item For $x\in H\bs M$ and $N$ the stratum through $x$, there exists a
neighborhood $U$ of $x$ in $H\bs M$ such that the isotropy equivalence class of
each element of the complement of $N$ in $U$ strictly dominates that of $x$.

\item If $M$ is compact, then the stratifications of $M$ and of $H\bs M$ are
finite.
\end{enumerate}
\end{prop}
 
 \begin{cor}\label{cor.orbstrata}  Let $\O$ be an orbifold.  Then the action of
$O(n)$ on the frame bundle $F(\O)$ gives rise to a (Whitney) stratification of
$\O$.  The strata are connected components of the isotropy equivalence classes
in $\O$.    The
set of regular points of $\O$ intersects each connected component $\O_0$ of $\O$
in a single stratum comprising an open dense submanifold of $\O_0$. 
 \end{cor}
 
 \begin{nota}\label{nota.ustrata} \text{ }
\begin{enumerate}
\item We will refer to the strata of $\O$ in Corollary \ref{cor.orbstrata} as
\emph{$\O$-strata}.
 
\item If $(\tu,G_U,\pi_U)$ is an orbifold chart on $\O$, then the action of
$G_U$
on $\tu$ gives rise to stratifications both of $\tu$ and of $U$ as in
Proposition \ref{prop.strat}.  We will refer to these as \emph{$\tu$-strata} and
\emph{$U$-strata}, respectively.
\end{enumerate}
 \end{nota}
 
 \begin{prop}\label{prop.ustrata}  Let $\O$ be a Riemannian orbifold and
$(\tu,G_U,\pi_U)$  be an orbifold chart.  Then:
 \begin{enumerate}
 \item
The $U$-strata are precisely the connected components of the intersections of
the $\O$-strata with $U$.  

\item Any two elements of the same $\tu$-stratum have the same stabilizers
 in $G_U$ (not just conjugate stabilizers).  

\item  If $H$ is a subgroup of $G_U$, then each connected component $W$ of the
fixed point set $\Fix(H)$ of $H$ in $\tu$ is a closed submanifold of $\tu$.  Any
$\tu$-stratum  that intersects $W$  nontrivially lies entirely in $W$.  Thus the
stratification of $\tu$ restricts to a stratification of $W$.
\end{enumerate}
\end{prop}

\begin{proof} (i) A consequence of Proposition \ref{prop.strat} is that each
$U$-stratum, respectively $\O$-stratum, is a connected component of the set of
all points in $U$, respectively $\O$, having $G_U$-isotropy, respectively
$O(n)$-isotropy, of a given order.  Since by \ref{localisot}, the  order of the
$G_U$-isotropy of each $x\in U$ is equal to the order of the $O(n)$-isotropy of
$x$, statement (i) follows.  

(ii) This follows because $G_U$ is discrete and the $\tu$-strata are
connected.

(iii)  The first statement is true for the fixed point set of any smooth proper
action by a compact group on a manifold \cite{DK}.  The second statement
follows from (ii). 
 \end{proof}

\begin{notarem}\label{nota.ists}  
Let $\O$ be a Riemannian orbifold and $(\tu,G_U,\pi_U)$  be an orbifold chart.
Let $\tn$ be a $\tu$-stratum in $\tu$.  By Proposition \ref{prop.ustrata}, all
the points in $\tn$ have the same isotropy group in $G_U$; we will refer to this
group as the isotropy group of $\tn$, denoted $\operatorname{{\Iso}}(\tn)$.

Given a $\tu$-stratum $\tn$, denote by $\ismax(\tn)$ the set of all $\gamma\in
{\Iso}(\tn)$ such that $\tn$ is open in the fixed point set $\Fix(\gamma)$ of
$\gamma$.

For $\gamma\in G_U$, Proposition \ref{prop.ustrata} tells us that each component
$W$ of the fixed point set $\Fix(\gamma)$ of $\gamma$ (equivalently, the fixed
point set of the cyclic group generated by $\gamma$) is a manifold stratified by
a collection of $\tu$-strata.   By Remark \ref{rem.strat}(i), the strata in $W$
of maximal dimension are open and their union has full measure in $W$.  In
particular, the union of those $\tu$-strata $\tn$ for which $\g\in\ismax(\tn)$
has full measure in $\Fix(\g)$.
 \end{notarem}

\begin{exa}\label{exa:square}
On $\R^2$, let $r_x$ and $r_y$ denote the reflections across the
$x$-axis and $y$-axis, respectively, and let $r_0$ denote the rotation through
angle $\pi$ about the origin $0$.  Then $G:=\{r_x,r_y,r_0,Id\}$ is a Klein four
group acting isometrically on $\R^2$.  The quotient of $\R^2$ by the semi-direct
product of $G$ with the lattice $\Z^2$ of translations is a closed orbifold
$\O$, whose underlying space is a square of side length $\frac{1}{2}$.  The
points on the boundary of the square are singular points (not boundary points)
of the orbifold, comprising eight strata:  each corner point forms a single
stratum with isotropy of order four, while each open edge forms a stratum with
isotropy of order two.  The strata of codimension one are called \emph{reflectors} or
\emph{mirrors}, and the
single point strata are called \emph{dihedral points} or \emph{corner reflectors}.  The
intersection $U$ of the square with a disk of radius less than $\frac{1}{2}$
centered at one of the corners is the image of an orbifold chart $(\tu,
G,\pi_U)$ where $\tu$ is a disk in $\R^2$ centered at the origin and $G$ is the
Klein-four group above.  The $\tu$-strata of this action consist of the single
point $0$ and 
the intersections of the disk $\tu$ with the positive and negative $x$-axis and
the positive and negative $y$-axis.  If $\tn$ is the intersection of $\tu$ with
one of the half axes, then ${\Iso}(\tn)$ consists of a
reflection and the identity, while $\ismax(\tn)$ contains only the reflection. 
For $\tn=\{0\}$, we have ${\Iso}(\tn)=G$, but $\ismax(\tn)=\{r_0\}$.
\end{exa}

%%%%%%%%%%%%%%%%%%%%%%%%%%%%%%%%%%%%%%%%%%%%%%%%%%%%%%%%%%%%%%%%%%%
%%%%%%%%%%%%%%%%%%%%%%%%%%%%%%%%%%%%%%%%%%%%%%%%%%%%%%%%%%%%%%%%%%%
\section{Construction of the heat expansion}\label{sec.heatkernel}
%%%%%%%%%%%%%%%%%%%%%%%%%%%%%%%%%%%%%%%%%%%%%%%%%%%%%%%%%%%%%%%%%%%
%%%%%%%%%%%%%%%%%%%%%%%%%%%%%%%%%%%%%%%%%%%%%%%%%%%%%%%%%%%%%%%%%%%

In this section, we address the heat kernel on closed Riemannian orbifolds. 

\begin{prop}\label{prop.spec}  Let $\O$ be a closed Riemannian orbifold. 
The Laplacian $\Delta$ of $\O$ has a discrete spectrum
$\lambda_1\leq\lambda_2\leq\dots$, with each eigenvalue having finite
multiplicity.  The normalized eigenfunctions $\varphi_j$ are $C^\infty$ and form
an orthonormal basis of $L^2(\O)$.  
\end{prop}

This proposition was proved in the case of $V$-manifolds (as defined in Remark
\ref{rem.embed}) by Y.-J. Chiang in \cite{Chi90}.  For an orbifold $\O$ which is
not a 
$V$-manifold, those strata of the singular set of codimension one are called
\emph{reflectors}.  
By doubling along all reflectors one obtains a $V$-manifold $X$ that doubly
covers $\O$.  
Thus $\O$ is the quotient of a $V$-manifold $X$ by a $\Z_2$ action. 
Since the eigenfunctions on $\O$ are then the $\Z_2$-invariant eigenfunctions on
$X$, Proposition \ref{prop.spec} follows immediately.

\begin{dfn}\label{def.kernel} Set $$\rp=[0,\infty)$$
and $$\rps=(0,\infty).$$  We say that  $K:\rps\times\O\times\O\to\R$ is a
\emph{fundamental solution} of the heat equation, or \emph{heat kernel}, if it
satisfies:
\begin{enumerate}

\item $K$ is $C^0$ in the three variables, $C^1$ in the
first, and $C^2$ in the second;

\item $\heat K(t,x,y)=0$ where $\Delta_x$ is the Laplacian with respect to the
second variable;

\item $\displaystyle \lim_\too\,K(t,x,\cdot)=\delta_x$ for all $x\in\O$.

\end{enumerate}
\end{dfn}

By the same argument as in the manifold case (see \cite{BGM}, III.E.2),
Proposition \ref{prop.spec} implies:

\begin{cor}\label{cor.unique}  
If a heat kernel exists, then it is unique and is given by 
$$K(t,x,y)=\sum_{j=1}^\infty\,e^{-\lambda_jt}\varphi_j(x)\varphi_j(y).$$
\end{cor}

  Chiang proved the existence of the heat kernel on a compact
$V$-manifold (from which existence on an arbitrary closed orbifold trivially
follows)  by proving the convergence of
$\sum_{j=1}^\infty\,e^{-\lambda_jt}\varphi_j(x)\varphi_j(y)$.  
She also showed that the heat kernel can be  approximated on good neighborhoods
by the Dirichlet heat kernel on the local manifold covering.  The existence also
follows from more general results on existence of the heat kernel for the basic
Laplacian on Riemannian foliations \cite{PR}.  
However, in
order to apply Donnelly's results on the heat trace for good orbifolds to obtain
the terms in the asymptotic expansion of arbitrary orbifolds in an applicable
form, we will not assume
the earlier existence results for the heat kernel or heat trace.  We instead
construct a
parametrix and then follow the standard construction of the heat kernel from the
parametrix as in \cite{BGM}.  Our construction of the parametrix and
consequently the heat kernel will use directly the local structure of orbifolds
as quotients of manifolds by finite group actions.

\begin{dfn}\label{def.param}
A \emph{parametrix} for the heat operator on $\O$ is a function
$H:\rps\times\O\times\O\to\R$ satisfying:
\begin{enumerate}
\item $H\in C^\infty(\rps\times\O\times\O)$;
\item $\heat H(t,x,y)$ extends to a function in $C^0(\rp\times\O\times\O)$;
\item $\displaystyle \lim_\too\,H(t,x,\cdot)=\delta_x$ for all $x\in\O$.
\end{enumerate}
\end{dfn}

Recall that the heat kernel on a closed $n$-dimensional Riemannian manifold
$M$ has an asymptotic expansion along the diagonal in $M\times M$ as $\too$ of
the form
\begin{equation}\label{eq.man.exp} 
K(t,x,x)\sim (4\pi
t)^{-\frac{n}{2}}(u_0(x,x)+tu_1(x,x)+t^2u_2(x,x)+\dots)
\end{equation}
where the $u_i$ are local Riemannian invariants defined in a neighborhood of the
diagonal in $M\times M$.  Letting $\zeta$ be a cut-off function that is
identically one near the diagonal, then for $m>\frac{n}{2}$, the function
\begin{equation}\label{eq.par}
K^{(m)}(t,x,y)= \zeta(x,y)(4\pi
t)^{-\frac{n}{2}}e^{-\frac{d(x,y)^2}{4t}}(u_0(x,y)+\dots + t^m
u_m(x,y))
\end{equation}
is a parametrix for the heat operator on $M$.

\begin{rem} \label{rem.loc} In what follows, we shall take a local covering
of our orbifold $\O$ by distinguished charts and piece together a parametrix for
the heat operator on $\O$ from the expressions in (\ref{eq.par}).  The key to
piecing together the parametrix on $\O$ is to note that while the parametrix
$K^{(m)}$ in (\ref{eq.par}) is globally defined on $M$, the three defining
conditions of a parametrix in Definition \ref{def.param} are satisfied locally
as well as globally by $K^{(m)}$.   Indeed the first two conditions are
trivially local.  The third condition is local in the following sense:  The
expression
$e^{-\frac{d(x,y)^2}{4t}}$ goes to zero uniformly as $\too$ when $d(x,y)$ is
bounded away from zero.  Thus for $f\in C^0(M)$, we have 
$$f(x)= \lim_\too \int_M\,K^{(m)}(t,x,y)f(y)dy=\lim_\too
\int_W\,K^{(m)}(t,x,y)f(y)dy,$$
where $W$ is any neighborhood of $x$.

We will also use the fact (see \cite{BGM}) that if $m>\frac{n}{2}+2l$, then
$\heat K^{(m)}(t,x,y)$  extends  to a function in $C^l(\rp\times M\times M)$. 
Moreover, the extension
is of class $C^{2l}$ in the last two variables.\end{rem}

\begin{nota}\label{nota.hm} Let $\O$ be an orbifold of dimension $n$.  Fix
$\e>0$ so that for each $x\in\O$,
there exists a distinguished coordinate chart of radius $\epsilon$ centered at
$x$.   Cover $\O$ with finitely many such charts
$(\twa, G_\a,\pi_\a)$, $\alpha=1,\dots, s$.  (Here we write $G_\a$ for
$G_{\wa}$ and $\pi_\a$ for $\pi_\wa$.)  Let
$\pa$ be the center of $\wa$ and $\tpa$ the center of $\twa$.   Let $\ua$,
respectively $\va$, be the geodesic ball of radius $\frac{\e}{4}$, respectively
$\frac{\e}{2}$, centered at $\pa$, and let $\tua$ and $\tva$ be the
corresponding
balls centered at $\tpa$ in $\twa$.  We may assume that the family of balls
$\{\ua\}_{1\le\a\le s}$ still covers $\O$.

For each $\a$ and each nonnegative integer $m$, we define  $\tham:
\rps\times\twa\times\twa$ by
$$\tham(t,\tx,\ty)=\sha\eha$$
where the $u_i$ are the invariants in \eqref{eq.man.exp}.
Since each $\g\in\ga$ is an isometry of $\twa$, we have
$u_i(\g\tx,\g\ty)=u_i(\tx,\ty)$ for all $\tx,\ty\in\twa$.  It follows that the
function
$$(t,\tx,\ty)\mapsto\sum_{\g\in\ga}\tham(t,\tx,\g\ty)$$
is $\ga$-invariant in both $\tx$ and $\ty$ and
thus descends to a well-defined function, which we denote by $\ham$,
on
$\rps\times\wa\times\wa$.

Let $\psia:\O\to\R$ be a $C^\infty$ cut-off function, which is identically one
on $\va$
and is supported in $\wa$.  Let $\{\etaa :\a=1,\dots,s\}$ be a partition of
unity on
$\O$ with the support of $\etaa$ contained in $\overline{\ua}$.   Define
$\hm:\rps\times\O\times\O\to\R$ by 
\begin{equation}\label{eq.hm}\hm(t,x,y)=\sum_{\a
=1}^s\,\psia(x)\etaa(y)\ham(t,x,y).\end{equation}  

\end{nota}

We will show that $\hm$ is a parametrix for the heat kernel on $\O$ when
$m>\frac{n}{2}$.

\begin{lem}\label{lem.cond1}  $\hm\in C^\infty(\rps\times\O\times\O)$.
\end{lem}  

Lemma \ref{lem.cond1} is immediate.

\begin{lem}\label{lem.cond2}  Let $l$ be a nonnegative integer.  Then 
\begin{enumerate}
\item $\heat\hm(t,x,y)$ extends to a function in
$C^l(\rp\times\O\times\O)$ if $m>\frac{n}{2}+2l$.  (It is moreover of class
$C^{2l}$ in the last two variables.)
\item For any given $T>0$ and for each $m>\frac{n}{2}$, there exists a constant
$A$ such that $|\heat\hm(t,x,y)|<At^{m-\frac{n}{2}}$ when $0<t <T$. 
\end{enumerate} 
\end{lem}

\begin{proof}   (i) Let $l\geq 0$ and suppose $m>\frac{n}{2}+2l$.  By Remark
\ref{rem.loc}, the function $\theat\tham(t,\tx,\ty)$ on
$\rps\times\twa\times\twa$ extends to $C^l(\rp\times\twa\times\twa)$.   (Here,
we are using the notation $\widetilde{\Delta}$ for the Laplacian on $\twa$ for
all choices of $\a$.)
   Thus the same is true for $\sum_{\g\in\ga}\theat\tham(t,\tx,\g\ty)$, and
hence
it follows that $\heat\ham(t,x,y)$ extends to 
   $C^l(\rp\times\wa\times\wa)$.   
   
   Now consider the function
$$f_\a(t,x,y):=\heat(\psia(x)\etaa(y)\ham(t,x,y)).$$  Noting that $\psia$ and
$\etaa$ are compactly supported inside $\wa$, we may view $f_\a$ as a function
on $\rps\times\O\times\O$ which is zero whenever $x$ or $y$ lies outside of
$\wa$. We show $f_\a$ extends to $C^l(\rp\times\O\times\O)$. Since $\psia\equiv
1$ on
$\va$, it follows immediately from the previous paragraph that
$f_\a$ extends to $C^l(\rp\times\va\times\O)$.  Moreover $f_\a\equiv 0$ on
$\rps\times\O\times(\O\setminus\ua)$ and so $f_\a$ also extends smoothly (to
zero) on $\rp\times \O\times(\O\setminus\ua)$.  Finally for $(x,y)\in
(\O\setminus\va)\times\overline{\ua}$, we have $d(x,y)\geq\frac{\epsilon}{4}$. 
Thus, as $\too$, $f_\a(t,x,y)$
and all its derivatives converge to zero uniformly for
$(x,y)\in(\O\setminus\va)\times\overline{\ua}$.  Thus $\hm$ extends to a
function in
$C^l(\rp\times\O\times\O)$.  The parenthetical statement in (1) similarly
follows from Remark \ref{rem.loc}.

(ii)  The $u_i$ are constructed so that 
$$\theat\tham(t,\tx,\ty)=\sha t^m\widetilde{\Delta}_\tx u_m(t,\tx,\ty)$$
(see \cite{BGM}).   Since the $u_i$ are $C^\infty$ functions, it follows that
there exists a constant $B_\a$ such that
$$\theat\tham(t,\tx,\ty)<B_\a t^{m-\frac{n}{2}}$$ on
$(0,T]\times\twa\times\twa$. 
Consequently, $$\heat\ham(t,x,y)<|\ga|B_\a t^{m-\frac{n}{2}}$$ on
$(0,T]\times\wa\times\wa$ and $$\heat
(\psia(x)\etaa(y)\ham(t,x,y))<|\ga|B_\a t^{m-\frac{n}{2}}$$ on
$(0,T]\times\va\times\O$, since $\psia\equiv 1$ on $\va$.  Once again, for $x$
outside of $\va$ and $y$ in the
support of $\etaa$, we have $d(x,y)\geq\frac{\epsilon}{4}$, and thus $\heat
(\psia(x)\etaa(y)\ham(t,x,y))$  can be bounded in terms of any power of $t$ on
$(0,T)\times(\O\setminus\va)\times \O$.  

Statement (ii) now follows from (\ref{eq.hm}).

\end{proof}

\begin{rem}\label{rem.bound}  When $m>\frac{n}{2}+2l$, one obtains bounds on
the partial derivatives of order at most $l$ of $\heat\hm(t,x,y)$ by an argument
analogous to that used in the proof of Lemma \ref{lem.cond2}(ii).
\end{rem}

\begin{lem}\label{lem.cond3}  Let $m$ be any nonnegative integer.  For $f\in
C^\infty(\O)$ and $x\in\O$, we have 
$$\lim_\too \int_\O \hm(t,x,y)f(y)dy=f(x).$$
I.e., $\lim_\too\hm(t,x,y)=\dx(y).$  Moreover, if $N$ is a topological space and
$f$ is a continuous function on $N\times\O$, then the convergence of \\ $\int_\O
\hm(t,x,y)f(p,y)dy$ to $f(p,x)$ is locally uniform on $N\times\O$.  
\end{lem}

\begin{proof}  Let $\tpsia$ and $\tetaa$ be the lifts of $\psia$ and $\etaa$ to
$\twa$.

Since $\supp(\etaa)\subset\overline{\ua}\subset\wa$, we have
\begin{equation}\label{eq.ham}\begin{aligned}
&\int_\O\psia(x)\etaa(y)\ham(t,x,y)f(p,y)\,dy\\&
=\frac{\psia(x)}{|G_\a|}\sum_{\g\in\ga}\int_\twa\,\tham(t,\tx,\g\ty)
\tetaa(\ty)\tfa(p,\ty)d\ty
\end{aligned}\end{equation}
where $\tfa$
  is a
lift of $f_{|N\times\wa}$ to
$N\times\twa$ and $\tx$ is an arbitrarily chosen point in the preimage of $x$
under the map $\twa\to\wa$.  
  
We change variables in each of the integrals in the right side of 
\eqref{eq.ham}, letting $\tul=\g(\ty)$.  Since
$\g$ is an isometry and since $\tetaa$ and $\tfa(p,\cdot)$
are $\g$-invariant, each integral in the summand is equal to 
\begin{equation}\label{eq.tham}\int_\twa\,\tham(t,\tx,\tul)
\tetaa(\tul)\tfa(p,\tul)\,d\tul.\end{equation}
As $\too$, the
integral (\ref{eq.tham}) above converges to
$\tetaa(\tx)\tfa(p,\tx)=\etaa(x)f(p,x)$ (see Remark \ref{rem.loc}).  Moreover,
this convergence is locally uniform on $N\times\twa$ (see \cite{BGM}).  Noting
that $\psia\equiv 1$ on the support of $\etaa$, it follows that  both sides of
\eqref{eq.ham} converge to $\etaa(x)f(p,x)$ as $\too$, and the
convergence is locally uniform on $N\times\wa$.    Since both sides of
\eqref{eq.ham} are identically zero when $x$ lies outside of
$\supp(\psia)\subset \wa$, we thus have locally uniform convergence to
$\etaa(x)f(p,x)$ on all of $N\times\O$.
  
Finally it follows from \eqref{eq.hm} that 
$$\lim_\too\int_\O\,\hm(t,x,y)f(p,y)\,dy=\sum_\a\,\etaa(x)f(p,x)=f(p,x)$$
locally uniformly.
\end{proof}

\begin{prop}\label{prop.param} 
$H^{(m)}$ is a parametrix for the heat
operator on $\O$ if $m>\frac{n}{2}$.
\end{prop}

\begin{proof}  
Immediate from Lemmas \ref{lem.cond1}, \ref{lem.cond2}(i), and
\ref{lem.cond3}.
\end{proof}

The construction of the heat kernel from the parametrix $\hm$ follows exactly as
in \cite{BGM}.  We give only a brief summary.

\begin{nota}\label{conv}  For $A,B\in C^0(\rp\times\O\times\O)$, define the
convolution $A*B\in
C^0(\rps\times\O\times\O)$ by 
$$A*B(t,x,z)=\int_0^t\,d\theta\int_\O\,A(t-\theta,x,y)B(\theta,y,z)dy.$$
Note that the convolution operator $*$  is associative.
\end{nota}

\begin{lem}\label{lem.qm}  Let $l$ be any nonnegative integer, and let
$m>\frac{n}{2}+2l$.  Define $F_m(t,x,y)
=\heat\hm(t,x,y)$.  (See Lemma \ref{lem.cond2} for regularity properties of
$F_m$.)  Then for each $T>0$, the series
$\sum_{j=1}^\infty\,(-1)^{j+1}F_m^{*j}(t,x,y)$
converges uniformly on $[0,T]\times\O\times\O$.  Let
$Q_m:\rp\times\O\times\O\to\R$ be the sum of this series.  Then $Q_m\in
C^l(\rp\times\O\times\O)$.  Moreover,  for any $T>0$, there exists a constant
$C$ such that  
$$|Q_m(t,x,y)|\leq Ct^{m-\frac{n}{2}}$$  on $[0,T]\times \O\times\O$.

\end{lem}

The proof of Lemma \ref{lem.qm} is identical to that of Lemma E.III.6 of
\cite{BGM} and uses only Lemma \ref{lem.cond2}(ii) and Remark
\ref{rem.bound}.

\begin{lem}\label{lem.conv}  Let $m>\frac{n}{2}$.  For $P\in
C^0(\rp\times\O\times\O)$, the function $\hm *P$, defined formally by the
expression in Notation \ref{conv}, exists and is in $C^0(\rps\times\O\times\O)$. 
Moreover, if $m>\frac{n}{2}+l$, then $\hm *P$ is of class $C^l$ in the second
variable.  For $m >\frac{n}{2}+2$,  $\heat(\hm*P(t,x,y))$ exists and equals
$(P+\hm *P)(t,x,y).$
\end{lem}

Again the proof is identical to that of Lemma E.III.7 of \cite{BGM} and is based
on Lemma \ref{lem.cond3}.  

Using Lemmas \ref{lem.qm} and \ref{lem.conv}, we obtain as in Proposition
E.III.8 of \cite{BGM} that:

\begin{prop}\label{prop.ker} Let $m>\frac{n}{2}+2$ and define $Q_m$ as in
Lemma
\ref{lem.qm}.   Then $K:=\hm -\hm *Q_m$ is a fundamental solution of the heat
equation on $\O$.
\end{prop}
  
  Note that uniqueness of the heat kernel implies that $\hm-\hm*Q_m$ is
independent of the choice of $m>\frac{n}{2}+2$.

\begin{nota}\label{nota.ha}  Let 
$$\tha(t,\tx,\ty)=\sum_{\g\in\ga}\shag\ehagi.$$
Observe that $\tha$ is $\ga$-invariant in both $\tx$ and $\ty$ and
thus descends to a well-defined function, which we denote by $\ha$,
on $\rps\times\wa\times\wa$.  
\end{nota}

\begin{thm}\label{thm.heattrace}  In the notation of Proposition \ref{prop.spec}
and \ref{nota.ha}, the trace of the heat kernel has an asymptotic expansion as
$\too$ given by 
$$\sum_{j=1}^\infty\,e^{-\lambda_jt}
\sim_{\too}\,\sum_{\a=1}^s\,\int_\O\,\etaa(x)\ha(t,x,x)\,dx.$$
\end{thm}

\begin{proof}  By Corollary \ref{cor.unique} and Proposition \ref{prop.ker}, we
have
  $$\sum_{j=1}^\infty\,e^{-\lambda_jt}=\int_\O\,(\hm-\hm*Q_m)(t,x,x)dx$$
  for $m>\frac{n}{2}+2$.  By Lemma \ref{lem.qm} and the fact that $(4\pi
t)^{\frac{n}{2}}\hm(t,x,x)$ is uniformly bounded for  $(t,x)\in(0,T]\times\O$
for any given $T>0$, it follows that 
  \begin{equation}\label{eq.int} \sum_{j=1}^\infty\,e^{-\lambda_jt}
  =\int_\O\,\hm(t,x,x)dx + O(t^{m-n})\end{equation}
  on any interval $(0,T]$.   (Aside:  When $\O$ is a manifold, then 
$\int_\O\hm(t,x,x)dx=(4\pi t)^{-\frac{n}{2}}(a_0+a_1t+\dots+a_mt^m)$.  For
general orbifolds, the arguments in the next section will show that
$\int_\O\hm(t,x,x)dx$ is of the form $(4\pi
t)^{-\frac{n}{2}}\sum_{j=0}^{2m}c_jt^{\frac{j}{2}}$.  Thus the error term can be
improved to $O(t^{m-\frac{n-1}{2}})$
since $\int_\O\,\hm(t,x,x)dx=\int_\O
H^{(m+n)}(t,x,x)dx +O(t^{ m-\frac{n-1}{2}})$.)
  
  Since $\psia$ is identically one on the support of $\etaa$, Notation
\ref{nota.hm} yields
 
\begin{equation}\label{eq.nopsi}\hm(t,x,x)=\sum_\a\,\etaa(x)\ham(t,x,x).
\end{equation}
  Substituting \eqref{eq.nopsi} into \eqref{eq.int}, we obtain
the theorem.
  \end{proof}

%%%%%%%%%%%%%%%%%%%%%%%%%%%%%%%%%%%%%%%%%%%%%%%%%%%%%%%%%%%%%%%%%%%
%%%%%%%%%%%%%%%%%%%%%%%%%%%%%%%%%%%%%%%%%%%%%%%%%%%%%%%%%%%%%%%%%%%
\section{Computation of the heat asymptotics}\label{sec.asymptotics}
%%%%%%%%%%%%%%%%%%%%%%%%%%%%%%%%%%%%%%%%%%%%%%%%%%%%%%%%%%%%%%%%%%%
%%%%%%%%%%%%%%%%%%%%%%%%%%%%%%%%%%%%%%%%%%%%%%%%%%%%%%%%%%%%%%%%%%%

\begin{notarem}\label{notarem.bk}  Let $\g$ be an isometry of a Riemannian
manifold $M$ and let $\om(\g)$ denote the set of components of the fixed point
set of $\g$.  Each element of $\om(\g)$ is a submanifold of $M$.  For each
non-negative integer $k$, Donnelly \cite{D1} defined a real-valued function, 
which we temporarily denote $b_k((M,\g),\cdot)$, on the fixed point set of $\g$. 
For each $W\in\om(\g)$, the restriction of $b_k((M,\g),\cdot)$ to $W$ is smooth. 
Two key properties of the $b_k$ are:

\begin{itemize}
\item (Locality)  For $a\in W$,   $b_k((M,\g),a)$ depends only  on the germs at
$a$ of
 the Riemannian metric of $M$ and of the isometry $\g$.  In particular, if $U$
is a $\g$-invariant neighborhood of $a$ in $M$, then
$b_k((M,\g),a)=b_k((U,\g),a).$  

\item (Universality)  If $M$ and $M'$ are Riemannian manifolds admitting
isometries $\g$ and $\g'$, respectively,  and if
$\sigma:M\to M'$ is an isometry satisfying $\sigma\circ \g=\g'\circ \sigma$, 
then
$b_k((M,\g),x)=b_k((M',\g'),\sigma(x))$ for all $x\in \operatorname{Fix}(\g)$.  

\end{itemize}

In view of the locality property, we will usually delete the explicit reference
to $M$ and rewrite these functions as $b_k(\g,\cdot)$, as they are written in
\cite{D1}.  

\end{notarem}

\begin{bk}\label{bk.comp}\cite{D1}

In the notation of \ref{notarem.bk}, let  $W\in\om(\g)$, and let $n=\dim(M)$ and
$m=\dim(W)$.  For $x\in W$, the orthogonal complement $T_x(W)^\perp$ of $T_x(W)$
in the tangent space $T_x(M)$ is invariant
under $\g_*$.  Define $A_\g(x)=\g_*:T_x(W)^\perp\to
T_x(W)^\perp$, and observe that $A_\g(x)$ is nonsingular.  Set
$$B_\g(x)=(I-A_\g(x))^{-1}.$$  
Donnelly 
showed that  $$b_k(\g,x)=|\det(B_\g(x))|\tb_k(\g,x),$$ where $\tb_k(\g,\cdot)$
is an
$O(m)\times O(n-m)$ universal invariant polynomial in the components of $B_\g$
and in the curvature tensor $R$ of $M$ and its covariant derivatives.

Explicit formulae for $b_0$ and $b_1$ are given in
\cite[Thm. 5.1]{D1} using the following indexing
conventions:  $1\leq \alpha,\beta\leq m$,  $m+1\leq i,j,k\leq n$ and $1\leq
a,b,c\leq n$.  At each point $x\in W$, choose an orthonormal basis
$\{e_1,\dots,e_n\}$ of $T_x(M)$ so that the first $m$ vectors are tangent to
$W$.  The sign convention on the curvature tensor $R$ of $M$ is chosen so that
$R_{abab}$  is the sectional curvature of the plane spanned by $e_a$ and
$e_b$.  Set
$$\tau=\sum_{a,b=1}^n\,R_{abab}$$
and $$\rho_{ab}=\sum_{c=1}^n\,R_{acbc}.$$
Thus $\tau$ is the scalar curvature and $\rho$ the Ricci tensor of $M$.  Then
\begin{equation}\label{eq.b0}b_0(\g,x)=|\det(B_\g(x))|\end{equation}
and, summing over repeated indices, 
\begin{equation}\label{eq.b1}\begin{aligned}
b_1(\g,x)=|\det(B_\g(x))|&(\frac{1}{6}\tau+\frac{1}{6}\rho_{kk}+\frac{1}{3}
R_{iksh}B_{ki}
B_{hs}\\+&\frac{1}{3}R_{ikth}B_{kt}B_{hi}-R_{kaha}B_{ks}B_{hs}).
\end{aligned}\end{equation}
\end{bk}

\begin{nota}\label{nota.orbbk}   Let $\O$ be an orbifold and let
$(\tu,G_U,\pi_U)$ be an orbifold chart.  

 In the notation of \ref{nota.ustrata} and \ref{nota.ists}, let $\tn$ be a
$\tu$-stratum and  let $\g\in\ismax(\tn)$.  Then $\tn$ is an open subset of a
component of $\operatorname{Fix}(\g)$ and thus by \ref{notarem.bk},
$b_k(\g,\cdot)$ ($=b_k((\tu,\g),\cdot)$) is smooth on $\tn$ for each nonnegative
integer $k$.  Define a function $b_k(\tn,\cdot)$ on $\tn$ by 
$$b_k(\tn,x)=\sum_{\g\in\ismax(\tn)}\,b_k(\g,x).$$

\end{nota}

\begin{lem}\label{lem.bk}  Let $\O$ be a Riemannian orbifold, let $N$ be an
$\O$-stratum and let $p\in N$.  Let $(\tu,G_U,\pi_U)$ and
$(\tu',G_{U'},\pi_{U'})$ be two orbifold charts with $p\in U\cap U'$.    Let
$\tp\in\tu$ and $\tp'\in\tu '$ with $\pi_U(\tp)=p=\pi_{U'}(\tp')$, and let
$\tn$, respectively $\tn'$, be the $\tu$-stratum through $\tp$, respectively
$\tu'$-stratum through $\tp'$.    Then for each $k$, we have
$b_k(\tn,\tp)=b_k(\tn',\tp')$.  
\end{lem}  

\begin{proof}  By Definition \ref{def.orb}, it suffices to consider the case
that one chart embeds in the other, say $\lambda: (\tu,G_U,\pi_U)\to
(\tu',G_{U'},\pi_{U'})$ is an isometric embedding with $\lambda(\tp)=\tp'$.  The
associated homomorphism $\tau:G_U\to G_{U'}$ carries ${\Iso}_{G_U}(\tp)$
isomorphically onto ${\Iso}_{G_{U'}}(\tp')$  and $\ismax(\tp)$ to
$\ismax(\tp')$.   The $\tu$-stratum $\tn$ is carried to an open subset of the
$\tu'$-stratum $\tn'$.  The lemma is thus an immediate consequence of  the
universality of the $b_k$, as discussed in \ref{notarem.bk}.
\end{proof}

\begin{dfn}\label{def.orbbk}   Let $\O$ be a Riemannian orbifold and let $N$ be
an $\O$-stratum.  

(i) For each non-negative integer $k$, define a real-valued function
$b_k(N,\cdot)$ by setting $b_k(N,p)=b_k(\tn,\tp)$ where $(\tu,G_U,\pi_U)$ is any
orbifold chart about $p$, $\tp\in\pi_U^{-1}(p)$ and $\tn$ is the $\tu$-stratum
through $\tp$.  By Lemma \ref{lem.bk}, the function $b_k(N,\cdot)$ is
well-defined.

(ii) The Riemannian metric on $\O$
induces a Riemannian metric, and thus a volume element, on the manifold $N$. 
Set
$$I_{N}:=(4\pi
t)^{-\dim(N)/2}\sum_{k=0}^\infty\,t^k\int_N\,b_k(N,x)d\vol_N(x)$$
where $d\vol_N$ is the Riemannian volume element. 

(iii) Also set
$$I_0=(4\pi t)^{-\dim(\O)/2}\sum_{k=0}^\infty\,a_k(\O)t^k$$
where the $a_k(\O)$ (which we will usually write simply as $a_k$) are the
familiar heat invariants.  More precisely, the invariants $u_i$ in
\eqref{eq.man.exp}, which are defined in terms of the curvature and its
covariant derivatives on any Riemannian manifold, also make sense on any
Riemannian orbifold.  The invariants $a_k(\O)$ are given by
$a_k=\int_\O\,u_k(x,x)d\vol_\O(x)$.  In particular, $a_0=\vol(\O)$,
$a_1=\frac{1}{6}\int_\O\tau(x)d\vol_\O(x)$, etc.   Note that if $\O$ is finitely
covered by a Riemannian manifold $M$, say $\O=G\bs M$, then
$a_k(\O)=\frac{1}{|G|}a_k(M)$.

\end{dfn}

\begin{thm}\label{thm.trace} Let $\O$ be a Riemannian orbifold and let
$\lambda_1\leq \lambda_2\leq\dots$ be the spectrum of the associated Laplacian
acting on smooth functions on $\O$.  The heat trace
$\sum_{j=1}^\infty\,e^{-\lambda_jt}$ of $\O$ is asymptotic as $\too$ to 
$$I_0+\sum_{N\in
S(\O)}\,\frac{I_N}{|\gn|}$$
where $S(\O)$ is the set of all $\O$-strata and where $|\gn|$ is the order of
the isotropy  at each $p\in N$ as defined in Remark \ref{localisot}.  This
asymptotic expansion is of
the form
\begin{equation}\label{eq.expform}(4\pi
t)^{-\dim(\O)/2}\sum_{j=0}^\infty\,c_jt^{\frac{j}{2}}\end{equation}
for some constants $c_j$.
\end{thm}

\begin{rem}\label{rem.goodcase}   Suppose $\O=G\bs M$ is a good closed
orbifold.  Note that $M$ may  be noncompact and $G$ may be an infinite group,
although the isotropy group at any point of $M$ must be a finite subgroup of
$G$.  In  this setting, Donnelly \cite{D2} proved the existence and
uniqueness of the heat kernel $K^M$ on $M$ and of an asymptotic expansion for
$K^M$.  He then obtained an asymptotic expansion for the heat trace on $\O$. 
Theorem \ref{thm.trace}, in the case of good orbifolds, organizes the
information in \cite{D2} in a way that clarifies the contribution of each
$\O$-stratum to the asymptotics.  

The expression for the heat asymptotics of good orbifolds in 
\cite{D2} differs from that in \eqref{eq.expform} in that the half
powers are
missing.  However, the absence of these powers is apparently a typographical
error
in transcribing a result of Donnelly's earlier paper \cite{D1}, stated below as 
Proposition \ref{prop.don}.

\end{rem}

K. Richardson \cite{Ric} obtained an asymptotic expansion for  the heat trace
associated with the basic Laplacian on a Riemannian foliation.  Also referring
to Donnelly's work on good orbifolds, he showed that the expansion is of the
form given in \eqref{eq.expform}.

The remainder of this section is devoted to the proof of Theorem
\ref{thm.trace}.

\begin{prop}\label{prop.don}\cite{D1}.  Let $M$ be a closed Riemannian manifold,
let $K(t,x,y)$ be the heat kernel of $M$, and let $\g$ be a nontrivial isometry
of $M$. 
Then, in the notation of \ref{notarem.bk},  $\int_M\,K(t,x,\g(x))d\vol_M(x)$ is
asymptotic as $t\to 0^+$ to
$$\sum_{W\in\om(\g)}\,(4\pi
t)^{-\frac{\dim(W)}{2}}\sum_{k=0}^\infty\,t^k\int_W\,b_k(\g,a)d\vol_W(a)$$
where $d\vol_W$ is the volume form on $W$ defined by the Riemannian metric
induced from $M$.
 \end{prop}
 
\begin{specpf}\label{specialpf}  We prove Theorem \ref{thm.trace} for
$\O=G\bs M$ a good closed orbifold with $G$ finite (and thus
$M$ compact).  In particular, $(M,G,\pi)$ is a global orbifold chart where
$\pi:M\to\O$ is the projection.  In this case, the theorem is an easy
consequence of Proposition
\ref{prop.don}.  Indeed, letting $K$ denote the heat kernel of $M$ and letting
$\pi:M\to\O$ be the projection, then the heat kernel $K^\O$ of $\O$ is given by
$$K^\O(t,x,y)=\sum_{\g\in G}\,K(t,\tx,\g(\ty))$$
where $\tx$, respectively $\ty$, are any elements of $\pi^{-1}(x)$, respectively
$\pi^{-1}(y)$.  Thus 
$$\int_\O K^\O(t,x,x) d\vol_\O(x)=\frac{1}{|G|}\sum_{\g\in
G}\,\int_M\,K(t,\tx,\g(\tx))d\vol_M(\tx),$$ so Proposition \ref{prop.don}
implies
that 
\begin{equation}\label{eq1}\begin{aligned}&\int_\O K^\O(t,x,x)
d\vol_\O(x)\sim_{\too}\frac{1}{|G|}\,\int_M\,K(t,\tx,\tx)d\vol_M(\tx)\\&+\frac{1
}{|G|}\sum_{1\neq \g\in G}\sum_{W\in\om(\g)}\,
(4\pi t)^{-\frac{\dim(W)}{2}}
\sum_{k=0}^\infty\,t^k\int_Wb_k(\g,a)d\vol_W(a).\end{aligned}\end{equation}

The first term on the right-hand-side of \eqref{eq1} is given by
\begin{equation}\label{eqfirst}
\frac{1}{|G|}\,\int_M\,K(t,\tx,\tx)d\vol_M(\tx)=
\frac{1}{|G|}(4\pi t)^{-\frac{\dim(M)}{2}}\sum_{k=0}^\infty\, a_k(M)t^k
=I_0.
\end{equation}
(See the final comment in Definition \ref{def.orbbk}(iii).)

Next let $1\neq \g\in G$, let $W\in\om(\g)$, and let $\tn$ be an $M$-stratum
contained in $W$.  Then either $\tn$ has measure zero in $W$ (in which case
$\g\notin\ismax(\tn)$) or else $\tn$ is open in $W$ and $\g\in \ismax(\tn)$. 
Thus by replacing the integral over $W$ with the integrals over the $M$-strata
that are open in $W$, reordering the summations in \eqref{eq1}, and
taking note of \eqref{eqfirst}, we obtain 
\begin{equation}\label{eqsecond}\int_\O K^\O(t,x,x)
d\vol_\O(x)\sim_{\too}I_0+\frac{1}{|G|}\sum_{\tn\in\widetilde{S}(M)}\widetilde{I
}_\tn\end{equation}
where $\widetilde{S}(M)$ denotes the set of all $M$-strata  and
where
$$\widetilde{I}_\tn=(4\pi
t)^{-\frac{\dim(\tn)}{2}}\sum_{k=0}^\infty\,t^k\int_\tn
b_k(\tn,a)d\vol_\tn(a).$$

Let $N$ be an $\O$-stratum.  Then $\pi^{-1}(N)$ is a union of finitely many
mutually isometric $M$-stratum $\tn_1,\dots, \tn_k$ and $\pi: \pi^{-1}(N)\to N$
is a covering map of degree $\frac{|G|}{|\Iso(N)|}$.    We have
$$\widetilde{I}_{\tn_1}+\dots + \widetilde{I}_{\tn_k}=\frac{|G|}{|\Iso(N)|}
I_N$$ and thus  Theorem
\ref{thm.trace}, in the case of orbifolds finitely covered by manifolds, follows
from \eqref{eqsecond}.
\end{specpf}

The proof in the general case will apply the argument in  \ref{specialpf} to
orbifold charts and then piece the computations together via a partition of
unity.  We first generalize Proposition \ref{prop.don} slightly.   The manifolds
in the two lemmas below do not have boundaries but could, for example, be
bounded domains in a larger manifold.

\begin{lem}\label{lem.temp} We use the notation of \ref{notarem.bk}. 
Let
$M$
be an $n$-dimensional
Riemannian manifold (without boundary) of finite volume and let $\g:M\to M$ be a
nontrivial isometry.  Assume that the distance $d(\tx,\g(\tx))$ remains bounded
away from zero off arbitrarily small tubular neighborhoods of the fixed point
set of
$\g$ and that each component of the fixed point set of $\g$ has finite volume. 
Then as $\too$,
$$\begin{aligned} \int_M\,\heatg\, d\tx\\\sim \sum_{W\in \om(\g)}(4\pi
t)^{-\dim(W)/2}\sum_{k=0}^\infty\,t^k\int_W\,b_k(\g,\tx)d\vol_W(\tx)
\end{aligned}$$
\end{lem}

This result is proven in Donnelly \cite[Thm. 4.1]{D1}, in case $M$ is
closed.   In that case, of course, the hypotheses on the distance function and
on the
fixed point set of $\g$ are automatic, and the lemma is a restatement of
Proposition \ref{prop.don}.  The proof goes through verbatim in the
more general setting of Lemma \ref{lem.temp}.  

\begin{lem}\label{lem.bketa}  With the notation and hypotheses of the previous
lemma,
let $\teta$ be a smooth bounded $\g$-invariant function on $M$.  Then there
exists a family of functions $c_k(\g,\teta,\cdot)$, $k=0,1,2,\dots$, defined on
the fixed point set of $\g$ and smooth on each component $W\in\om(\g)$, such
that as $\too$,
$$\begin{aligned} \int_M\,\teta(\tx)\heatg\, d\tx\\\sim \sum_{W\in \om(\g)}(4\pi
t)^{-\dim(W)/2}\sum_{k=0}^\infty\,t^k\int_W\,c_k(\g,\teta, \tx)d\vol_W(\tx).
\end{aligned}$$
Moreover,    $\ct_k(\g,\teta, \cdot)$
satisfies the following:
\begin{enumerate}
\item (Locality) $\ct_k(\g,\teta,\tx)$ depends only on the germs of $\g$,
$\teta$, and the
Riemannian metric of $M$ at  $\tx\in W$;
\item $\ct_k(\g,\teta, \cdot)$ is zero off $supp(\teta)\cap W$;  
\item the dependence of $\ct_k(\g,\teta,\cdot)$ on $\teta$ is linear;
\item $\ct_k(\g,1,\cdot)=b_k(\g,\cdot)$ where $1$ denotes the constant function
$\teta\equiv 1$;
\item (Universality) if $M'$ is another Riemannian manifold, $\g'$ is an
isometry of $M'$ and
$\sigma:M\to M'$ is an isometry satisfying $\sigma\circ \g=\g'\circ\sigma$, then
$\ct_k(\g',\teta\circ\sigma^{-1},\sigma(\tx))=\ct_k(\g,\teta,\tx)$ for all $\tx$
in the fixed point set of $\g$.
\end{enumerate}
\end{lem}

The proof requires only minor changes in the proof of Theorem 4.1 of
\cite{D1}. 
In the proof of that theorem, the functions $b_k(\g,\cdot)$ are expressed as
linear
combinations of certain derivatives of explicitly defined functions $h_j$,
$j=0,\dots,k$.  To obtain the functions $\ct_k(\g,\teta,\cdot)$, one replaces
the
functions
$h_j$ by the functions $\teta h_j$.  

\begin{notarem}\label{nota.ck}  Let $\O$ be a closed orbifold and consider the
charts $\tua,\tva,\twa$ and partition of 
unity $\{\etaa\}$ given in Notation \ref{nota.hm}.    Let $\tva$ play the role
of
$M$ in Lemmas \ref{lem.temp} and \ref{lem.bketa}, and  let
$\tetaa=\etaa\circ\pi_\a$ 
play the role of $\teta$.   Since $\tva$ has compact closure
inside the
larger Riemannian manifold $\twa$ on which $G_\a$ acts by isometries and since
the fixed point set of $\g$ in $\twa$ is connected (in fact, it is the union of a
collection of geodesics radiating from the center point $\tp_\a$), one easily
verifies for each $\g\in G_\a$ that  the hypothesis concerning the distance
function in the
two lemmas holds.  
\begin{enumerate}
\item For each $\tva$-stratum $\tn$, define a smooth function
$c_{k,\a}(\tn,\cdot)$ on $\tn$ by
$$c_{k,\a}(\tn,\tx)=\sum_{\g\in\ismax(\tn)}\,c_k(\g,\tetaa,\tx).$$

\item Let $N$ be an $\O$-stratum.  For each non-negative integer $k$ and each
$\alpha=1,\dots, s$, define a continuous (in fact, smooth) function
$c_{k,\a}(N,\cdot)$ on $N$ as
follows:  First for $x\in N\cap\va$, set  $c_{k,\a}(N,x)=c_k(\tn,\tx)$ where
$\tx$ is any element of $\pi_\a^{-1}(x)$ and $\tn$  is the $\tva$-stratum
through $\tx$.   By an argument analogous to that of Lemma \ref{lem.bk}. this
definition is independent of the choice of $\tx$ in $\pi_\a^{-1}(x)$.   Since 
$\tetaa$ is supported in $\overline{\tua}$, Lemma
\ref{lem.bketa}(ii) implies that $c_{k,\a}(N,\cdot)$ is zero off $N\cap
\overline{\ua}$ and thus extends to a continuous function on $N$ which is zero
off $N\cap \ua$.

\item Set
$$I_{N,\a}:=(4\pi
t)^{-\dim(N)/2}\sum_{k=0}^\infty\,t^k\int_N\,c_{k,\a}(N,x)d\vol_N(x).$$
\end{enumerate}
\end{notarem}

\begin{lem}\label{lem.ck}  Let $\O$ be a closed Riemannian orbifold and let $N$
be an $\O$-stratum.  Then for each non-negative integer $k$, we have
$$\sum_{\a=1}^s\,c_{k,\a}(N,\cdot)=b_k(N,\cdot)$$
and
$$\sum_{\a=1}^s\,I_{N,\a}=I_N.$$
\end{lem}

\begin{proof}  Let $x\in N$, and let $\a_1,\cdots,\a_r$ be those
$\a\in\{1,\dots,s\}$ for which $x\in V_{\a_i}$.  Then we can find a coordinate
chart $(\tu,G_U,\pi_U)$ such that $U\subset V_{\a_1}\cap\dots\cap V_{\a_r}$ and
such that the chart $(\tu,G_U,\pi_U)$ embeds in each of the charts
$(\tva,G_\a,\pi_\a)$.  Let $\lambda_i:\tu\to\tv_{\a_i}$ be the embedding.  Let
$\tx\in\pi_U^{-1}(x)$, let $\tn$ be the $\tu$-stratum through $\tx$, and let
$\tx_i=\lambda_i(\tx)$.  As in the proof of Lemma \ref{lem.bk}, $\lambda_i(\tn)$
is an open subset of the $\tv_{\a_i}$-stratum $\tn_i$ through $\tx_i$.    Using
the universality property (v) of Lemma \ref{lem.bketa}, Notation \ref{nota.ck},
and an argument analogous to that of Lemma \ref{lem.bk}, we see that
$$c_{k,\a_i}(N,x):=c_{k,\a_i}(\tn_i,\tx_{i})=\sum_{\g\in\ismax(\tn)}\,c_k(\g,\eta_{\a_i} \circ\pi_U,\tx).$$
Thus, since $c_{k,\a}(N,x)=0$ when $\a$ is not one of $\a_1,\dots,\a_r$, we have
\begin{equation}\label{eq.a}\sum_{\a=1}^s\,c_{k,\a}(N,x)=\sum_{\g\in\ismax(\tn)}
\,\sum_{i=1}^r\,c_k(\g,\eta_{\a_i}\circ\pi_U,\tx).\end{equation}
  From properties (iii) and (iv)  of Lemma \ref{lem.bketa} and the fact that
$\displaystyle\sum_{i=1}^r\eta_{\a_i}\equiv 1$ on $U$,  we see that  
\begin{equation}\label{eq.b}\sum_{i=1}^r\,c_k(\g,\eta_{\a_i}\circ\pi_U,\tx)=b_k(
\g,\tx)\end{equation}
on $\tn$.    By Definition \ref{def.orbbk},
  $$\sum_{\g\in\ismax(\tn)}b_k(\g,\tx)=b_k(N, x)$$
  and thus the first equation in the lemma follows from \eqref{eq.a}
and \eqref{eq.b}.  The second equation is then immediate.
 \end{proof}

We now prove Theorem \ref{thm.trace}.

\begin{proof} Let $n=\dim(\O)$.  By Theorem \ref{thm.heattrace} and the fact
that
the support of $\etaa$ is contained in $\va$, we have
\begin{equation}\label{eq.0} \sum_{j=1}^\infty\,e^{-\lambda_jt}\sim_\too
\sum_{\a=1}^s\,\int_\va\etaa(x)\ha(t,x,x)d\vol(\O)\end{equation}
where $\ha(t,x,x)$ is defined in Notation  \ref{nota.ha}.  By Notation
\ref{nota.ha},
\begin{equation}\label{eq.mess}\begin{aligned}&\sum_{\a=1}^s\,\int_\va\,\phia(x)
\ha(t,x,x)\,d
\vol_\O(x)=\\
&\sum_{\a=1}^s\frac{1}{|\ga|}\int_\tva\tetaa(\tx)(4\pi
t)^{-\frac{n}{2}}\ehaxi\,d\vol_\tva(\tx) \\
&+\sum_{\a=1}^s\frac{1}{|\ga|}\sum_{1\neq \g\in G_\a}\int_\tva\tetaa(\tx)\shagx
(u_0(\tx,\g(\tx))\\&\hspace{2in}+tu_1(\tx,\g(\tx))+\dots)\,d\vol_\tva(\tx)
\end{aligned}\end{equation}
where $\tetaa=\etaa\circ\pi_\a$. 

Consider the first sum on the right-hand-side of \eqref{eq.mess}. 
Since $\etaa$ is supported in $\va$, we
have
\begin{equation}\label{eq.term1}\begin{aligned}&\sum_{\a=1}^s\,\frac{1}{|\ga|}(4
\pi
t)^{-\frac{n}{2}}\int_\tva\,\tetaa(\tx)\ehaxi\,d\vol_\tva(\tx)\\
=&\sum_{\a=1}^s(4\pi
t)^{-\frac{n}{2}}\int_\O\,\etaa(x)(u_0(x,x)+tu_1(x,x)+\dots)d\vol_\O(x)\\
&=(4\pi t)^{-\frac{n}{2}}\int_\O\,(u_0(x,x)+tu_1(x,x)+\dots)d\vol_\O(x) =I_0.
\end{aligned}\end{equation}

Next by Lemma \ref{lem.bketa} and the remarks in \ref{nota.ck}, we have for each
$1\neq\g\in\ga$,
\begin{equation}\label{eq.term2}\begin{aligned}&\frac{1}{|\ga|}\sum_{1\neq\g\in\
ga}\int_\tva\,\tetaa(\tx)\shagx
(u_0(\tx,\g(\tx))+\\&\hspace{1.5in}tu_1(\tx,\g(\tx))+\dots)\,d\vol_\tva(\tx)\\
\end{aligned}\end{equation}
$$\sim\frac{1}{|\ga|}\sum_{1\neq\g\in\ga}\sum_{W\in \om(\g)}(4\pi
t)^{-\dim(W)/2}\sum_{k=0}^\infty\,t^k\int_W\,\ct_k(\g,\tetaa,\tx)d\vol_W(\tx).$$

By \ref{nota.ck} and an argument identical to that in  \ref{specialpf}, the
right-hand-side of \eqref{eq.term2} is
equal to 
$$\sum_{N\in S(\O)}\,\frac{1}{|\Iso(N)|}I_{N,\a}.$$
Consequently, Lemma \ref{lem.ck} and \eqref{eq.term2} imply that
the second sum in the right-hand-side of \eqref{eq.mess} is equal to 
$$\sum_{N\in S(\O)}\,\frac{1}{|\Iso(N)|}I_{N}.$$
Thus in view of \eqref{eq.0}, \eqref{eq.mess},
and \eqref{eq.term1}, the theorem is proved.

\end{proof}

%-------------------------------------------
%        Begin Applications 
%-------------------------------------------

%%%%%%%%%%%%%%%%%%%%%%%%%%%%%%%%%%%%%%%
%%%%%%%%%%%%%%%%%%%%%%%%%%%%%%%%%%%%%%%
\section{Applications}\label{sec.apps}
%%%%%%%%%%%%%%%%%%%%%%%%%%%%%%%%%%%%%%%
%%%%%%%%%%%%%%%%%%%%%%%%%%%%%%%%%%%%%%%

\begin{thm} \label{4.1}  Let $\O$ be a Riemannian orbifold with singularities. 
If $\O$ is even dimensional (respectively, odd dimensional) and some
$\O$-stratum of the singular set is odd dimensional (respectively,
even dimensional),
then $\O$ cannot be isospectral to a Riemannian  manifold. 
\end{thm}

\begin{proof}   
It is clear from (\ref{eq.b0}) that if $N$ is any $\O$-stratum of
the singular
set, then the function $b_0(\g,\cdot)$ is strictly positive on $N$ for each
$\g\in\gns$.  Thus, in the two cases, the fact that $\O$ is an orbifold can be
gleaned from the presence of half-integer powers, respectively integer powers,
of $t$ in the asymptotic expansion in Theorem~\ref{thm.trace}.  
\end{proof}

\begin{rem}
In \cite{GR03} this theorem was stated for good orbifolds, but here we also
include bad orbifolds.
\end{rem}

We now restrict our attention to closed two-dimensional orbifolds 
(2-orbifolds).  The singularities which may occur in 2-orbifolds are 
cone points, dihedral corner reflectors, and mirror reflectors.  
Recall that dihedral corner reflectors and mirror reflectors both appeared in 
Example \ref{exa:square}.  A cone point $p$ of order $n$ is an isolated
singularity; 
an orbifold chart for a neighborhood of $p$ is $(\mathbb{D}^2, \mathbb{Z}_n,
\pi)$ 
where $\mathbb{D}^2$ is an open 2-disk in $\R^2$ and $\mathbb{Z}_n$ 
is the cyclic group of order $n$.  The Euler characteristic 
of a 2-orbifold is 2 minus the sum of the related values: each cone point of
order $n$ has value $\frac{n-1}{n}$; each dihedral corner reflector has value
$\frac{n-1}{2n}$; each handle has value 2; each cross-cap has value 1; and each
mirror reflector has value 1.  Every good 2-orbifold admits a
(metrically) spherical, Euclidean or hyperbolic structure depending on whether
the Euler characteristic is positive, zero  or negative, respectively
\cite{Thu80}.  In addition, all
bad 2-orbifolds have positive Euler characteristic.

%-------------------------------------------
%        4.2 Cone b
%-------------------------------------------
\begin{exa}\label{4.2}   Let $\O$ be a 2-orbifold and let
$p$ be a cone point of order $m$.  If $N=\{p\}$, then $\gn$ is
a cyclic group of order
$m$ and $\gns$ contains all of the nontrivial elements.  Letting $\g$ be the
generator, then for $j=1,\dots,m-1$ we have that
$$A_{\g^j}=\g^j_*=\begin{bmatrix} \cos(\frac{2j\pi}{m})&-\sin(\frac{2j\pi}{m})\\
\sin(\frac{2j\pi}{m})& \cos(\frac{2j\pi}{m})\end{bmatrix},$$
where $A_{\g^j}$ is as defined in 4.2.
Thus
$$b_0(\g^j)=|det((I-A_{\g^j})^{-1})|=\frac{1}{2-2\cos(\frac{2j\pi}{m})}=\frac{1}
{4\sin^2(
\frac{j\pi }{m})}.$$
(We are writing $b_0(\g^j)$ rather than using the function notation
$b_0(\g^j,\cdot)$
since $N$ consists of a single point.)
\end{exa}

%-------------------------------------------
%        4.3 Residue Lemma
%-------------------------------------------

\begin{lem}\label{4.3}
$$\sum_{j=1}^{m-1}\,\frac{1}{\sin^2(\frac{j\pi}{m})}=\frac{m^2-1}{3}.$$
\end{lem}

\begin{proof}  A well-known formula (see, for example, \cite[Example 7.9.1]{Pennisi}), proven by the calculus of residues, states that
$$\frac{\pi^2}{\sin^2(\pi z)}=\sum_{k=-\infty}^\infty\,\frac{1}{(k-z)^2}.$$
Thus
$$\sum_{j=1}^{m-1}\,\frac{1}{\sin^2(\frac{j\pi}{m})}=\frac{m^2}{\pi^2}\,\sum
_{j=1}^{m-1}\,\sum_{k=-\infty}^\infty\,\frac{1}{(mk-j)^2}.$$

Since
$$\sum_{j=1}^{m-1}\,\sum_{k=-\infty}^\infty\,\frac{1}{(mk-j)^2}=2\,\sum_{n
=1}^\infty\,\frac{1}{n^2}-2\,\sum_{n=1}^\infty\,\frac{1}{m^2n^2}$$
and $\sum_{n=1}^\infty\,\frac{1}{n^2}=\frac{\pi^2}{6}$, the lemma follows.
\end{proof}

%-------------------------------------------
%        4.4 I_N Computation for Cone Point
%-------------------------------------------

\begin{prop}\label{4.4}  Let $\O$ be a 2-orbifold, let $p$ be a cone point
in
$\O$ of order $m$ and let $N=\{p\}$.  Then in the notation of Theorem
~\ref{thm.trace}, we have
$$I_N=\frac{m^2-1}{12}+O(t).$$
\end{prop}

\begin{proof}  By \ref{def.orbbk}(ii) and \ref{4.2},  
$$I_N=\sum_{j=1}^{m-1}\,\frac{1}{4\sin^2(\frac{j\pi}{m})} +O(t).$$
Thus Proposition \ref{4.4}  follows from Lemma \ref{4.3}.
\end{proof}

\begin{exa}\label{2orbterms}Calculating heat invariants for 2-orbifolds.

%-------------------------------------------
%        4.5 Degree zero term for cone points
%-------------------------------------------

 {\bf Degree zero term for orientable 2-orbifolds}.
An orientable 2-orbifold
$\O$ can have only isolated singularities,  i.e., cone points.  Suppose $\O$ has
$k$ cone points of orders $m_1,\dots,m_k$.  In the notation of
\ref{def.orbbk}(iii), 
$$I_0=\frac{1}{4\pi}(a_0 t^{-1} +a_1 +O(t)).$$  Thus by Theorem~\ref{thm.trace}
and Proposition \ref{4.4}, the term of degree zero  in the asymptotic expansion
in Theorem~\ref{thm.trace} is given by
$$\frac{a_1}{4\pi} +\sum_{i=1}^{k}\frac{1}{m_i}\frac{m_i^2-1}{12}.$$
By the Gauss-Bonnet Theorem (valid also for orbifolds; see \cite{Sat57,Thu80}), we have 
$$a_1=\frac{2\pi}{3}\chi(\O).$$  
Hence the degree zero term is

\begin{equation} 
\label{eqn:deg0} \frac{\chi(\O)}{6}+\sum_{i=1}^{k}\frac{m_i^2-1}{12m_i}. 
\end{equation}

%-------------------------------------------
%        Degree zero term for nonorientable 2-orbifolds
%-------------------------------------------

  {\bf Degree zero term for nonorientable 2-orbifolds}.
For a 2-orbifold $\O$, the dimension zero singular locus is the only portion
which contributes to the degree zero term of the asymptotic expansion in 
Theorem~\ref{thm.trace}, aside
 from the $\frac{a_1}{4\pi} = \frac{\chi({\O})}{6}$ component.
A nonorientable 2-orbifold $\O$ can have cone points and/or 
dihedral corner reflector points that contribute in the following ways.
As in the computation of the degree zero term for orientable 2-orbifolds, 
here a simple cone point of order
$m$ contibutes $\frac{1}{m}\frac{m^2-1}{12}$.  
Let $N=\{p\}$, where $p$ is a corner reflector point created by  
a rotation of order $n$ and a reflection.  Then $|\gn| = 2n$.  By
Notation \ref{nota.ists}(ii), $\gns$ contains only the nontrivial elements of
the
rotation group, since the reflection fixes one-dimensional strata of
the mirror locus, a higher dimensional stratum of the singular set.  Hence the
computations in Example~\ref{4.2} and Proposition~\ref{4.4} remain the same, but
the difference in
$|\gn| $ is seen as an extra $\frac{1}{2}$ factor in Theorem~\ref{thm.trace}. 
The
point contributes
 $\frac{1}{2n}\frac{n^2-1}{12}$ to the degree zero term.
 
 Thus, for a 2-orbifold $\O$ with cone points $p_1, \ldots, p_k$ of orders
$m_1,\ldots, m_k$, 
and dihedral corner reflector points 
$q_1, \ldots, q_r$ of orders $n_1, \ldots,n_r$,
the term of degree 0 in the asymptotic expansion of the heat trace is
\begin{equation} \label{eqn:deg0rev} \frac{\chi({\O})}{6} + \sum_{i=1}^k
\frac{1}{m_i} \frac{m_i^2-1 }{12 } 
+ \sum_{j=1}^r \frac{1}{2n_j} \frac{n_j^2-1}{12}.
\end{equation}

%-------------------------------------------
%        4.5 Degree one term
%-------------------------------------------

 {\bf Degree one term for 2-orbifolds}.
Aside from $a_2$, only dimension zero strata of the singular set contribute to
the $t$ term of the asymptotic expansion in Theorem~\ref{thm.trace}.  We first
note that the last term in (\ref{eq.b1}) is zero; this follows from the
summation convention  and the fact that our singular set is zero-dimensional. 
Using symmetry properties of the curvature, we can further simplify
(\ref{eq.b1}) as 
\begin{equation*}
b_1(\gamma^j) = \frac{1}{4 \sin ^2 (\frac{j\pi}{m})}\left( \frac{\tau}{6} +
\frac{\rho_{kk}}{6} + \frac{2}{3} (R_{1212}(B_{21}^2 + B_{12}^2 - B_{12}B_{21} -
B_{22}B_{11}))\right),
\end{equation*}
where $R_{1212}$ is evaluated in the local covering manifold.
By the definitions of scalar and Ricci curvature, the preceding equation becomes
\begin{equation*}
b_1(\gamma^j) = \frac{R_{1212}(1 + B_{21}^2 + B_{12}^2 - B_{12}B_{21} -
B_{22}B_{11})}{6 \sin ^2 (\frac{j\pi}{m})}. 
\end{equation*}
In general, straightforward calculations show that 
$$
B_{\gamma^j} = (I-A_{\gamma^j})^{-1} = 
\begin{bmatrix}
\frac{1}{2} & -\frac{\sin (\frac{2j\pi}{m})}{2-2\cos (\frac{2j\pi}{m})} \\
\frac{\sin (\frac{2j\pi}{m})}{2-2\cos (\frac{2j\pi}{m})} & \frac{1}{2}
\end{bmatrix}
$$
which implies
\begin{equation*}\label{eqn:b_1final}
b_1(\gamma^j) = \frac{R_{1212}}{8 \sin ^4(\frac{j\pi}{m})},
\end{equation*}
for $j=1, \ldots, m-1$.

Thus, for a 2-orbifold $\O$ with cone points $p_1, \ldots, p_k$ of orders $m_1,
\ldots, m_k$, and dihedral corner reflector points $q_1, \ldots, q_r$ of orders
$n_1, \ldots, n_r$,
the coefficient of the term of degree 1 in the asymptotic expansion of the heat
trace is
\begin{equation}\label{eqn:degree1}
 \frac{a_2}{4\pi} + \sum_{i=1}^k \frac{1}{m_i} \left( \sum_{j=1}^{m_i-1}
\frac{R_{1212}}{8 \sin ^4 (\frac{j \pi}{m_i})} \right)
+ \sum_{i=1}^r \frac{1}{2n_i} \left( \sum_{j=1}^{n_i-1} \frac{R_{1212}}{8 \sin
^4 (\frac{j \pi}{n_i})} \right).
\end{equation}
Recall that $a_2(\O) = \frac{1}{360} \int_\O (2|R|^2 - 2 |\rho|^2 + 5 \tau^2)
d\vol_\O(g)$, where $R$ is the curvature, $\rho$ is the Ricci curvature, and
$\tau$ is the scalar curvature of $\O$ (e.g. \cite{BGM}).  

We can further simplify (\ref{eqn:degree1}) by making the substitution 
\begin{equation*}
\sum_{j=1}^{m_i-1} \frac{1}{\sin ^4 (\frac{j \pi}{m_i})} = 
\frac{m_i^4 +10m_i^2-11}{45}.
\end{equation*}
See \cite{BY02} and the references therein or \cite{Chencsc4,SIAMcsc4} for
evaluations of this and similar finite trigonometric sums.

 With this substitution, (\ref{eqn:degree1}) becomes
\begin{equation}\label{eqn:alpha2}
 \frac{a_2}{4\pi} + \sum_{i=1}^k \frac{R_{1212}(m_i^4 + 10m_i^2 - 11)}{360m_i}+
 \sum_{i=1}^r \frac{R_{1212}(n_i^4 + 10n_i^2 - 11)}{720n_i}.
\end{equation}
\end{exa}

%-------------------------------------------
%        Degree 1/sqrt(t)  Donnelly left these off
% in one paper as a typo only
%-------------------------------------------

{\bf Degree $-\frac{1}{2}$ term for 2-orbifolds}.
The only $\O$-strata that contribute to the degree $\frac{1}{\sqrt{t}}$ term are
those of codimension one in $\O$.  To obtain these $\O$-strata, remove any
dihedral
points from the mirror locus and then take the connected components of the
remaining set.  Let $x \in N$, an $\O$-stratum of the mirror locus, and note
that
$\g \in \gns$ must act as a reflection.  

To compute $b_0(\g,x)=|det((I-A_\g)^{-1})|$, notice that on the
normal bundle to  $N$,
$\g_* = [-1]$, and so, 
$$b_0(\g,x)=|det((I-A_\g)^{-1})| = | [2]^{-1}| = \frac{1}{2}.$$
Applying \ref{def.orbbk}(ii),
$$I_N= (4\pi t)^{-1/2}\sum_{\g\in\gns}\sum_{k=0}^\infty\,t^k\int_N
\frac{1}{2}d\vol_N(x) +O(t)$$
$$=\frac{length(N)}{4 \sqrt{\pi}} \frac{1}{\sqrt{t}} + O(\sqrt{t}).
$$
We sum over all $\O$-strata of the mirror locus to obtain the coefficient of the
$\frac{1}{\sqrt{t}}$
term in Theorem~\ref{thm.trace}:
\begin{equation} \label{eqn:deg-1/2} \sum_N\,\frac{I_N}{|\gn|} = 
\sum_N\,\frac{1}{2} \frac{length(N)}{4 \sqrt{\pi}} =
\frac{length(MirrorLocus(\O))}{8 \sqrt{\pi}}.\end{equation}

%-------------------------------------------
%        4.5 Degree sqrt(t) term
%-------------------------------------------

 {\bf Degree $\frac{1}{2}$ term for 2-orbifolds}.
For a 2-orbifold $\O$, the dimension one singular locus gives the sole
contribution to the $\sqrt{t}$ term of the asymptotic expansion in
Theorem~\ref{thm.trace}.

We have
\begin{eqnarray*}
b_1(\g,x) & = &
|det(B_\g(x))|(\frac{\tau}{6}+\frac{\rho_{kk}}{6}+
\frac{1}{3}R_{iksh}B_{ki} B_{hs} \\ 
& & \mbox{} 
+\frac{1}{3}R_{ikth}B_{kt}B_{hi}-R_{k\a h\a}B_{ks}B_{hs})(x)
\end{eqnarray*}
where $B_{ij}$ denotes the $i,j$ entry of $B_\g(x)$, $\tau$ is the scalar
curvature of $M$ at $x$ and $\rho$ is the Ricci curvature. 

In the case of a 2-orbifold with a point $x$ in its mirror locus (but not a
dihedral corner reflector point), 
the matrix $B_\g(x)$ is one-dimensional, and thus the
third and fourth terms in the sum vanish.  Since $B_\g(x)=\frac{1}{2}$, the last
term
is $-\frac{1}{4}R_{1212}$.  We also have $\rho_{kk}=R_{1212}$, while
$\tau=2R_{1212}$.  Thus
$b_1(\g,x)=\frac{1}{2}(\frac{1}{4}R_{1212})(x)=\frac{1}{8}R_{1212}(x)=
\frac{1}{16}\tau(x)$.  Thus for $N$ a stratum in the mirror locus, the term of
degree
$\frac{1}{2}$ in $I_N$ is given by 
$$\frac{\sqrt{t}}{\sqrt{4\pi}}\int_N
\frac{1}{16}\tau(x)d\vol_N(x)=
\frac{\sqrt{t}}{32\sqrt{\pi}}\int_N\tau(x)d\vol_N(x).$$
The coefficient of the degree $\frac{1}{2}$ term in the asymptotic expansion for
the heat trace on $\O$ is thus
\begin{equation}\label{eqn:deg1/2}
\frac{1}{64\sqrt{\pi}}\int_{MirrorLocus(\O)}\tau
\end{equation}
where the scalar curvature is the scalar curvature of $\O$ computed at points in
the mirror locus and the integral is with respect to the induced Riemannian
metric on the one-dimensional mirror locus.  See Table \ref{table:2-orbs} for the asymptotic expansions of the heat kernel for orbifolds $\O$ with $\chi(\O) \geq 0$.

 %-------------------------------------------
%        Table of Spherical and Euclidean 2-orbifold Expansions
%-------------------------------------------

\begin{table*}[htpb]\label{table:exp}
\begin{minipage}{\textwidth}
 \caption{2-orbifold expansions with $\chi(\O) \geq 0$} 
\begin{tabular}{ l @{\quad} l @{\quad}l @{\quad}l }
\boldmath
$\O \mbox{ \textbf{with} } \chi(\O) > 0$ 
\unboldmath 
 \footnotetext{
Here $m,n \geq 1$, $V=\frac{\vol(\O)}{4\pi}$ and
$ML=\frac{length(MirrorLocus(\O))}{8\sqrt{\pi}}$.  Note that *torus is an
annulus 
with two mirror reflector edges and *Klein bottle is a M\"{o}bius band with one 
mirror reflector edge.}
 &\textbf{Asymptotic Expansion}&\\
 $\O(m)$  & $V \frac{1}{t}+  \frac{1}{12}(2+m+ \frac{1}{m}) +O(t)$\\
 $\O(*m)$ & $V \frac{1}{t} + ML\frac{1}{\sqrt{t}}  + \frac{1}{24}(2+m+
\frac{1}{m}) +O(\sqrt{t})$\\
 $\O(m,n)$ &$V \frac{1}{t} +  \frac{1}{12}(m+n+ \frac{1}{m} +
\frac{1}{n})+O(t)$\\
 $\O(*m,n)$ &$V \frac{1}{t} + ML\frac{1}{\sqrt{t}} +\frac{1}{24}(m+n
+\frac{1}{m}+\frac{1}{m})+O(\sqrt{t})$\\
 $\O(m\times)$  &$V \frac{1}{t}+\frac{1}{12}(m+\frac{1}{m})+O(t)$\\
 $\O(m*)$  & $V \frac{1}{t} + ML\frac{1}{\sqrt{t}}
+\frac{1}{12}(m+\frac{1}{m})+O(\sqrt{t})$\\
 $\O(2,2,m)$  &$V \frac{1}{t}+\frac{1}{12}(3+m+\frac{1}{m}) +O(t)$\\
 $\O(*2,2,m), \O(2,{*}m)$ & $V \frac{1}{t} + ML\frac{1}{\sqrt{t}}
+\frac{1}{24}(3+m+\frac{1}{m}) +O(\sqrt{t})$ \\
  $\O(2,3,3)$ & $V \frac{1}{t}+\frac{43}{72} +O(t)$\\
 $\O(*2,3,3),\O(3,{*}2)$  & $V \frac{1}{t} + ML\frac{1}{\sqrt{t}} 
+\frac{43}{144}+O(\sqrt{t})$\\
 $ \O(2,3,4)$  &$V \frac{1}{t}+\frac{97}{144}  +O(t)$ \\
 $\O(*2,3,4)$  &$V \frac{1}{t} + ML\frac{1}{\sqrt{t}} +\frac{97}{288}
+O(\sqrt{t})$\\
 $\O(2,3,5)$ & $V \frac{1}{t}+\frac{271}{360} +O(t)$\\
 $\O(*2,3,5)$ &$V \frac{1}{t} + ML\frac{1}{\sqrt{t}}  +\frac{271}{720}
+O(\sqrt{t})$\\
&&\\ 
\boldmath
 $\O \mbox{ \textbf{with} } \chi(\O) = 0$
\unboldmath 
 &\textbf{Asymptotic Expansion}&\\
 torus, Klein bottle  & $V \frac{1}{t}+ O(t)$\\
 *torus, *Klein bottle &$V \frac{1}{t} + ML\frac{1}{\sqrt{t}} +O(\sqrt{t})$\\
 $\O(2,2,2,2)$  & $V \frac{1}{t}+\frac{1}{2} +O(t)$\\
$\O({*}2,2,2,2),\O(2,{*}2,2),$ & \\
\, \, $\O(2,2{*})$ &  $V \frac{1}{t} + ML\frac{1}{\sqrt{t}} 
+\frac{1}{4}+O(\sqrt{t})$\\
  $\O(2,2\times)$ &  $V \frac{1}{t} +\frac{1}{4}+O(t)$\\ 
$\O(2,4,4)$  & $V \frac{1}{t}+\frac{3}{4} +O(t)$\\
  $\O(*2,4,4),\O(4,{*}2)$ &  $V \frac{1}{t} + ML\frac{1}{\sqrt{t}}
+\frac{3}{8}+O(\sqrt{t})$\\
 $\O(3,3,3)$  & $V \frac{1}{t}+\frac{2}{3} +O(t)$\\
  $\O(*3,3,3),\O(3,{*}3)$ &  $V \frac{1}{t} + ML\frac{1}{\sqrt{t}} 
+\frac{1}{3}+O(\sqrt{t})$\\
  $\O(2,3,6)$  & $V \frac{1}{t}+\frac{5}{6} +O(t)$\\
  $\O(*2,3,6)$ &  $V \frac{1}{t} + ML\frac{1}{\sqrt{t}}
+\frac{5}{12}+O(\sqrt{t})$\\
\end{tabular}
\label{table:2-orbs}
 \end{minipage}
\end{table*}

Our notation for orbifolds is adapted from 
Conway's convention \cite{Con92}, with commas added for readability. 
Namely $\O(a,b,*c,d)$ denotes an orbifold with 
simple cone points of orders
$a$ and $b$, and dihedral corner reflectors of orders $2c$ and $2d$.  In
addition, 
$\O(n\times)$ is a disk with a simple interior cone point and the edge
identified via the antipodal map, while $\O(n*)$ is a disk with a simple
interior cone point and a mirror edge corresponding to a reflection.  A detailed
explanation of the orbifold notation can be found in \cite{Con92} and compared
with pictures  in \cite[pages 80--90]{Mon85}.    
A proof that these are
all of the
closed 2-orbifolds with $\chi(\O) \geq 0$ and that the only bad 
2-orbifolds are
$\O(m), \O(*m)$, $\O(m,n)$ and $\O(*m,n)$ (when $m>1$ and $m \neq n$),  can be
found in \cite{Thu80}. 

Let $\O$ be an orientable 2-orbifold with $k$ cone points of orders $m_1,
\ldots, m_k$, denoted $\O(m_1, \ldots, m_k)$,
and consider the quantity $c$ defined as $12$ times the degree zero term:
\begin{equation} \label{eqn:defofc}
c=2\chi(\O)+\sum_{i=1}^{k}\left(m_i-\frac{1}{m_i}\right). 
\end{equation}
This quantity is a spectral invariant; note that it depends only on the
topology,
not on the Riemannian metric.   For $\O(m_1, \ldots, m_k)$, we denote by $c(m_1,
\ldots, m_k)$ the associated spectral invariant.  We now investigate classes of
orientable 2-orbifolds for which $c$ is a complete topological invariant.  Although Theorem \ref{4.6} is a special case of Theorem \ref{thm:classC}, we begin with the more restricted class in order to give the reader intuition for the proof techniques used.

%-------------------------------------------
%        4.6 Football and Teardrop Theorem
%-------------------------------------------

\begin{thm}\label{4.6}  Within the class of all footballs (good or bad) and all
teardrops, the spectral invariant $c$ is a 
complete topological invariant. I.e.,
$c$ determines whether the orbifold is a football or teardrop and determines the
orders of the cone points. 
\end{thm}

\begin{proof} Denote by $\O(m)$ the teardrop with cone point of order $m$ and by
$\O(r,s)$ the football with cone points of orders $r$ and $s$. Let $c(m)$ and
$c(r,s)$ denote the invariant defined in \eqref{eqn:defofc} in the two
cases. Then $\O(m)$ has Euler characteristic
$1+\frac{1}{m}$ and thus 
$$c(m)=2+m+ \frac{1}{m}.$$  The football $\O(r,s)$ has Euler characteristic
$\frac{1}{r}+\frac{1}{s}$, so 
$$c(r,s)=r+s+ \frac{1}{r}+\frac{1}{s}.$$
The invariant is an integer only in the case of $\O(2,2)$, so the football
$\O(2,2)$ is spectrally distinguishable from the other footballs and teardrops. 
Thus for the remainder of the proof, all footballs will be assumed to have at
least one cone point of order strictly greater than 2.  

For teardrops, $c(m)$ trivially determines $m$.  We  next claim that footballs
are distinguishable from teardrops.  Indeed, suppose that $c(m)=c(r,s)$.  Then
$m+2=r+s$ and $\frac{1}{m}=\frac{1}{r}+\frac{1}{s}$.  The latter equation
implies that $m< min(r,s)$.  Since also $r,s\geq 2$, we have $2+m<r+s$, a
contradiction, thus proving the claim.

It remains only to show that for footballs,
$c(r,s)$ determines $r$ and $s$.  From $c(r,s)$, one can read off the quantities
$r+s$ and $\frac{1}{r}+\frac{1}{s}=\frac{r+s}{rs}$.  Hence $c(r,s)$ determines
both $r+s$ and $rs$, thus also $|r-s|$, since $(r-s)^2=(r+s)^2-4rs$.  Hence
$(r,s)$ is determined up to order, completing the proof.
\end{proof}

%-------------------------------------------
%     4.7 More topological results
%-------------------------------------------

%-------------------------------------------
%        Triangular Pillow Orbifold Table 
%-------------------------------------------

\begin{table}\begin{center}
\caption{Triangular Pillow Orbifolds}
\begin{tabular}{  l @{\quad} lc l cl @{\quad}l }
& $\O$  &$\chi(\O)$ &$c(\O)$ \\
$\chi(\O)>0$& $\O(2,2,2)$  & $\frac{1}{2}$ &$5\frac{1}{2}$\\
& $\O(2,2,m)$  & $\frac{1}{m}$ &$3+m+\frac{1}{m}$\\
& $\O(2,3,3)$ & $\frac{1}{6}$ &$7 \frac{1}{6}$\\
& $\O(2,3,4)$ & $\frac{1}{12}$  &$8 \frac{1}{12}$ \\
& $\O(2,3,5)$ & $\frac{1}{30}$ &$ 9 \frac{1}{30} $\\ 
$\chi(\O)=0$& $\O(3,3,3)$  & $ 0$ &$ 8$ \\
& $\O(2,4,4)$  & $ 0$ &$9$\\
& $\O(2,3,6)$  & $ 0$ &$ 10$ \\
$\chi(\O)<0$& $\O(3,3,4)$ & $-\frac{1}{12}$ &$8\frac{11}{12}$\\
& $\O(3,4,4)$ &$-\frac{1}{6}$ &$9\frac{5}{6}$\\
& $\O(3,3,5)$ &$-\frac{2}{15}$ & $9\frac{13}{15}$\\
& $\O(2,4,5)$ &$-\frac{1}{20}$ & $9\frac{19}{20}$\\
& \vdots & $-1<\chi(\O) < 0$ &$c(\O) > 10$ \\
\end{tabular}
\label{table:pillows}
\end{center}
\end{table}

\begin{thm}\label{thm:classC}
Let $C$ be the class consisting of all closed orientable 2-orbifolds with
$\chi(\O) \geq 0$.
The spectral invariant $c$ is a complete topological invariant within $C$ and
moreover, it distinguishes the elements of $C$ from smooth oriented closed
surfaces.
\end{thm}

\begin{proof}
We first consider the 2-orbifolds for which $c$ is an integer.  Note that among
teardrops, footballs, and triangular pillows, the only integer values are
$c(2,2)=5, c(2,3,6)=10, c(2,4,4)=9$, and $c(3,3,3)=8$ (cf. Tables
~\ref{table:2-orbs} and \ref{table:pillows}).  In addition, $c(2,2,2,2)= 6,
c(S^2) = 4$, and $c(T^2) = 0$.  Let $S_g$ be a Riemann surface of genus $g \geq
2$.  Then we also have $c(S_g)=4-4g$.  It is clear that the values of $c$ are
distinct in each case, so that the spectrum distinguishes these 2-orbifolds.  

For the rest of the proof, it suffices to consider orbifolds in $C$ with
$\chi(\O) > 0$ 
since these include all orbifolds within $C$ for which $c$ is not an integer.
Table \ref{table:pillows} lists the values of $c$ for these triangular pillows. 
Setting $c(2,3,3) = c(2,2,m)$ and solving for $m$ gives $m=\frac{25 \pm
\sqrt{481}}{12}$, which contradicts the assumption that $m$ is an integer. 
Similar calculations for $c(2,3,4)$ 
and $c(2,3,5)$
show that the spectrum distinguishes among triangular pillows with $\chi(\O)
>0$. 

By Theorem \ref{4.6}, $c$ distinguishes among teardrops and footballs.
We next show that $c$ distinguishes both teardrops and footballs from
triangular pillows with $\chi(\O) >0$.  We have $c(m) = 2+m+ \frac{1}{m}$, and
$c(p,q,r) = -2 + p + q + r + \frac{1}{p} + \frac{1}{q} + \frac{1}{r}$; setting
the integer and fractional parts of these equations equal gives
\begin{eqnarray*}
2 + m & = & -1 + p + q + r \\
\frac{1}{m} & = & \frac{1}{p} + \frac{1}{q} + \frac{1}{r} -1.
\end{eqnarray*}
Solving the first equation for $m$ and plugging the result into the second
equation yields
$$
0 = pr(p+r-3) + pq(p+q-3) + qr(q+r-3) + pqr(5-p-q-r).
$$ 
None of the possible triples $(p,q,r)$ satisfy this equation, showing that
teardrops are distinguished from  triangular pillows with $\chi(\O) >0$.  The
argument that $c$ distinguishes good footballs from these triangular pillows is
analogous.  

To see that $c$ distinguishes triangular pillows with $\chi(\O) >0$ from bad
footballs, we compare the respective integer and fractional parts of $c$ in each
case.  For example, comparing $c(2,3,3)$ and $c(r,s), r \neq s$ gives
\begin{eqnarray*}
7 & = &  r + s \\
\frac{1}{6} & = & \frac{1}{r} + \frac{1}{s},
\end{eqnarray*}
which implies $s^2 - 7s + 42 =0$; thus $s = \frac{7 \pm \sqrt{-119}}{2}$,
contradicting $s$ being an integer.  The calculations are similar for
$\O(2,3,4)$ and $\O(2,3,5)$, while for $\O(2,2,m)$ we have
\begin{eqnarray*}
3 + m & = & r + s\\
\frac{1}{m} & = & \frac{1}{r} + \frac{1}{s},
\end{eqnarray*} 
which implies $r(r-3) + s(s-3) + rs = 0$.  We can assume that one of $r$ and $s$
is strictly greater than 2, say $r$.  Thus $r-3 \geq 0$ and $s-3 \geq -1$, which
implies 
$$
r(r-3) + s(s-3) + rs \geq s(r-1) >0.
$$
Hence triangular pillows with $\chi(\O) > 0$ are distinguished from bad
footballs.
\end{proof}

%--------------------------------------------
%       4.8 Results requiring metric
%---------------------------------------------

\begin{rem}  
Notably absent from the class $C$ are triangular pillows with  $\chi(\O) < 0$. 
The invariant $c$ does not seem sufficiently strong to distinguish among these
triangular pillows.  However, as a special case of a result in \cite{EBDStroh},
the spectrum does determine the orders of the cone points in such a 2-orbifold,
provided that it is endowed with a metric of constant curvature $-1$.  On the
other hand, to distinguish, say, triangular pillows  with $\chi(\O) < 0$ from
triangular pillows with  $\chi(\O) > 0$, we do not need such a metric
assumption.

If $\O(p,q,r)$ is any triangular pillow with  $\chi(\O) < 0$ and
$\spec(\O(p,q,r)) = \spec(\O(2,2,m))$,  then setting the integer and
fractional parts of the respective values of $c$ equal, we have
\begin{eqnarray*}
m+3 & = & p+q+r-2 \\
\frac{1}{m} & = & \frac{1}{p} + \frac{1}{q} + \frac{1}{r}.
\end{eqnarray*}
Solving the first equation for $m$ and plugging the resulting value into the
second equation yields
\begin{equation}\label{eqn:sph_hyp}
2pqr + pr(p+r-5) + pq(p+q-5) + qr(q+r-5) = 0.
\end{equation}
Note that for a  triangular pillow with  $\chi(\O) < 0$, we must have
$\frac{1}{p} + \frac{1}{q} + \frac{1}{r} <1$, which implies that the sum of any
two of $p,q,r$ is at least 5.  Thus each term on the right-hand side of
(\ref{eqn:sph_hyp}) is nonnegative.  This contradiction then implies that
$\O(2,2,m)$ cannot be isospectral to such a triangular pillow.  For the
remaining triples where  $\chi(\O) > 0$, we set the integer and fractional parts
of the respective value of $c$ equal to those of a triangular pillow with 
$\chi(\O) < 0$, and note that there are no  $\chi(\O) < 0$ triples satisfying
these equations.  Thus $c$ distinguishes between triangular pillows with 
$\chi(\O) < 0$ and  $\chi(\O) > 0$.  A similar argument shows that $c$
distinguishes  triangular pillows with  $\chi(\O) < 0$ from teardrops; it is
clear that $c$ also distinguishes  triangular pillows with  $\chi(\O) < 0$ from
the smooth surfaces and the remaining elements of $C$, with the exception of
footballs.  In this last case, it seems that metric assumptions are again
necessary.
\end{rem}

\begin{figure}
\center
\includegraphics[width=3in]{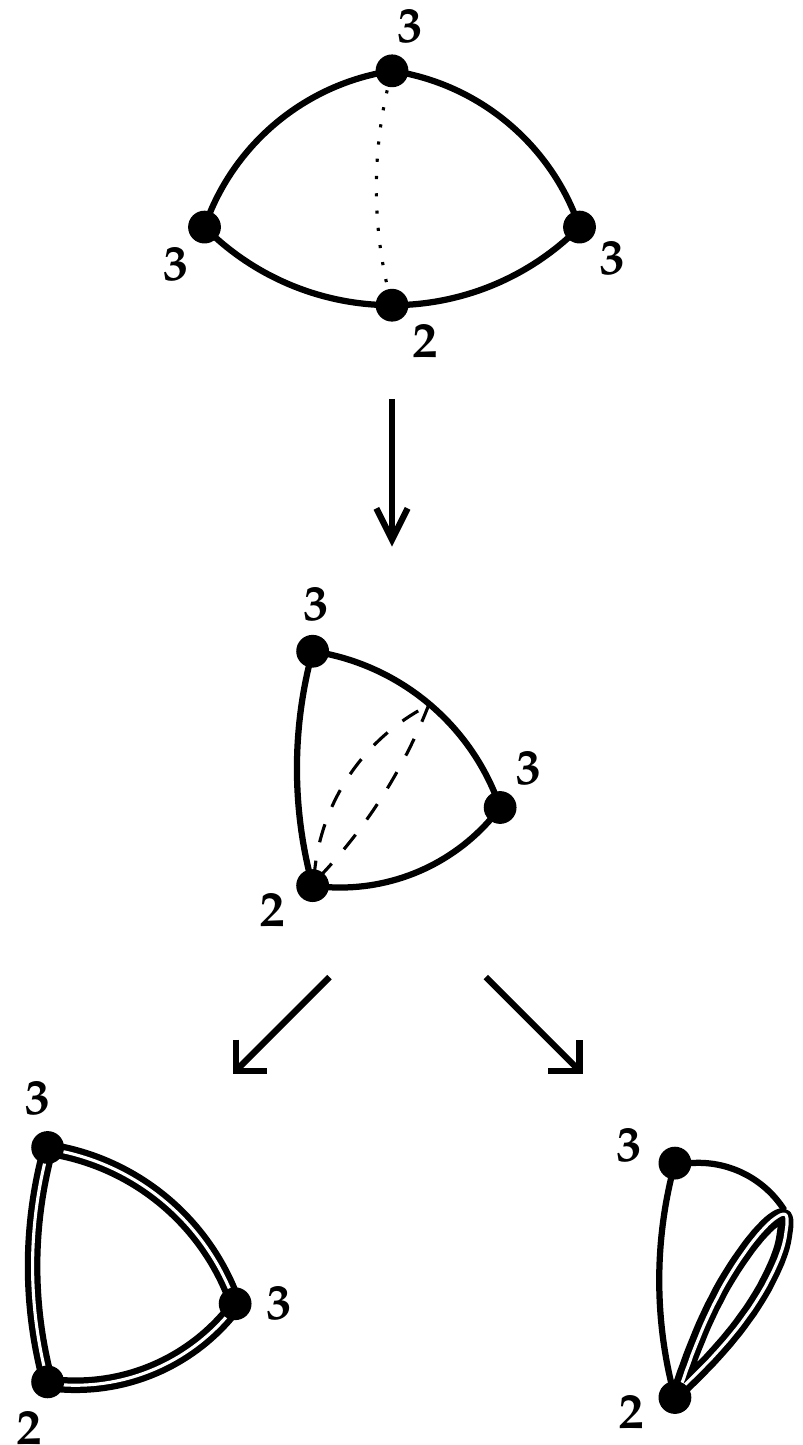}
\caption{The uppermost object is a fundamental domain for $\O (2,3,3)$, with
its quotient $\O(2,3,3)$ below.  Here the vertices are cone points labelled with
their orders, and the double edges represent reflector edges.  On the bottom
left we show $\O(*2,3,3)$, which is obtained by reflecting $\O(2,3,3)$ in the
plane of the paper, and on the right is $\O(3,*2)$, obtained by reflecting
$\O(2,3,3)$ in the plane containing the dashed loop.}
\label{fig:233}
\end{figure}

\begin{rem}
Other difficulties arise during the consideration of the expanded class which
includes nonorientable orbifolds.  Metric assumptions are necessary in numerous
cases, such as distinguishing nonorientable orbifolds with the same orientable
double cover (e.g. $\O(*2,3,3)$ and $\O(3,{*}2)$, see Table~\ref{table:2-orbs}
and Figure~\ref{fig:233}). 
\end{rem}

We now examine classes that include nonorientable orbifolds.

\begin{prop}  Within the class of all closed 2-orbifolds with $\chi(\O) \geq 0$,
the spectrum distinguishes whether the orbifold has zero or positive Euler
characteristic.
\end{prop}
\begin{proof}
Note that for 2-orbifolds $\O$ with $\chi(\O) = 0$, $c$ is either an integer or
equal to $4.5$.  Thus $c$ distinguishes all but the following cases:  $S^2$ from
the orbifolds $\O(*3,3,3)$ and $\O(3,{*}3)$ (with $c=4$), 
the good football $\O(2,2)$ from  $\O(*2,3,6)$ (with $c=5$), and
the bad teardrop $\O(2)$ from the orbifolds $\O(*2,4,4)$ and $\O(4,{*}2)$ (with
$c=4.5$).  The lack of a mirror locus in the $\chi(\O) > 0$ cases and the
presence of a mirror locus in the corresponding $\chi(\O) = 0$ orbifolds can be
gleaned from the degree $-\frac{1}{2}$ term, and so they are distinguished.
\end{proof}

\begin{prop}\label{hearingcurvaturesigncontinued}

Within the class of closed 2-orbifolds of constant nonzero curvature $R$ or
$-R$ the spectrum determines the sign of the curvature, i.e. whether the
orbifold is spherical or hyperbolic.

\end{prop}

\begin{proof}
Assume that $\O$ has cone points 
$p_1, \ldots, p_k$ of orders $m_1, \ldots, m_k$, 
and/or dihedral corner reflector points $q_1, \ldots, q_r$ of orders 
$n_1, \ldots, n_r$.
Now look at the coefficient of the $t$ term in the expansion, as in \eqref{eqn:alpha2},
which reduces to 
$$\frac{a_2}{4\pi} \pm R \left(\sum_{i=1}^k \frac{(m_i^2 + 11)
(m_i^2-1)}{360m_i}+
 \sum_{i=1}^r \frac{(n_i^2 + 11)(n_i^2-1)}{720n_i}\right)$$ in the presence of
constant curvature.  
The $\frac{1}{t}$ term in the expansion tells us $\vol(\O)$ and we know the size
of the curvature, so we know the $a_2$ component and may subtract it off. 
Notice that the summands are nonnegative, and hence we can read the sign of the
curvature unless there were no cone points and no dihedral corner reflector
points.  In this case, examine the degree zero term, which now has no point
contributions and reduces to 
$\frac{a_1}{4\pi} = \frac{\chi({\O})}{6}$, and we can read off the Euler
characteristic.
\end{proof}

\begin{rem}   In the case of closed 2-orbifolds with a
nontrivial mirror locus, \eqref{eqn:deg1/2} enables us to make a
stronger statement.  In particular, among such orbifolds that are endowed with
a metric of strictly positive, strictly negative, or zero curvature, 
the spectrum determines the sign of the curvature.  This class includes the 
bad orbifolds $\O(*m)$ and $\O(*p,q)$ with $p \neq q$ since they admit a metric
of strictly
positive (but variable) curvature.
\end{rem}

\begin{prop}\label{spherical}
Within the class of  spherical 2-orbifolds of constant curvature $R>0$ the
spectrum determines the orbifold.
\end{prop}

\begin{proof}
Notice that the metric  requirement eliminates teardrops and bad footballs and
their quotients from this class.
In  Table~\ref{table:2-orbs}, 
$c$ distinguishes among the remaining spherical orbifolds with the exception
that $c$ is unable to distinguish between orbifolds that are nonorientable but
have the same orientable double cover: $\O(*m,m), \O(m\times),$ and $\O(m*)$,
with double cover
$\O(m,m)$; $\O(*2,2,m)$ and $\O(2,*m)$, with double cover $\O(2,2,m)$; and
$\O(*2,3,3)$ and $\O(3,*2)$, with double cover $\O(3,3,2)$ (see Figure
\ref{fig:233}).  
However, $c$ is able to distinguish each nonorientable class from the
remaining orbifolds.

Consider such a class of nonorientable spherical orbifolds with a common
orientable double cover.  The  coefficient of the degree $-\frac{1}{2}$ term, as given in \eqref{eqn:deg-1/2}, distinguishes nonorientable
orbifolds with mirror loci from those without, i.e. in this class it
distinguishes orbifolds with only crosscaps from those with mirror loci.
In the presence of constant curvature, \eqref{eqn:deg-1/2} also
distinguishes among the remaining spherical cases:  among $\O(*m,m)$ and
$\O(m*)$, the length of the mirror locus of $\O(*m,m)$ is larger in the constant
curvature metric; among $\O(*2,2,m)$ and $\O(2,{*}m)$, the length of the mirror
locus of $\O(*2,2,m)$ is larger; and among $\O(*2,3,3)$ and $\O(3,{*}2)$, the
length of the mirror locus of $\O(*2,3,3)$ is larger (see Figure \ref{fig:233}).
\end{proof}

\begin{rem}  Notice that metric assumptions are only needed to distinguish
within each nonorientable class.  We cannot make a similar statement for flat
2-orbifolds.  For example, it is possible to endow $\O(2,*2,2)$ and $\O(2,2*)$
with a metric of zero curvature so that they have the same area and also have
mirror loci of the same length.  They cannot be distinguished by the asymptotic
expansion of the heat trace.
\end{rem}

%-------------------------------------------
%        Begin References
%-------------------------------------------

\end{document}